\documentclass[11pt, leqno]{article}

\usepackage{euscript,amstext,amsmath,amsthm,a4,amssymb}
\def\bC{{\rm \bf C}}

\def\bQ{{\rm \bf Q}}
\def\bR{{\rm \bf R}}

\def\bZ{{\rm \bf Z}}
\def\ch{{\rm  ch}}

\def\rk{{\rm  rk}}
\def\End{{\rm  End}}
\def\ind{{\rm  Ind}}
\def\dip{\displaystyle}

\newcommand{\Ker}{\rm Ker}
\newcommand{\tr}{{\rm Tr}}
\newcommand{\rL}{{\rm L}}
\newcommand{\sig}{{\rm sig}}

\newcommand{\var}{\varepsilon}
\newcommand{\wi}{ \widetilde}

\newcommand{\comment}[1]{}

\newtheorem{thm}{Theorem}[section]
\newtheorem{cj}[thm]{Conjecture}
\newtheorem{prop}[thm]{Proposition}
 \newtheorem{cor}[thm]{Corollary}
\newtheorem{lemma}[thm]{Lemma}
\theoremstyle{remark}
\newtheorem{Rem}[thm]{Remark}
\theoremstyle{definition}
\newtheorem{defn}[thm]{Definition}

\title{Eta-invariants, torsion forms\\ and flat vector bundles}
\author{Xiaonan Ma\footnote{Centre de Math\'ematiques, UMR 7640 du CNRS,
\'Ecole Polytechnique, 91128 Palaiseau Cedex, France.
(ma@math.polytechnique.fr)}\ \ and\ \ Weiping Zhang\footnote{Nankai Institute of Mathematics, Nankai
University, Tianjin 300071, P.R. China.
(weiping@nankai.edu.cn)}}
\date{}

\begin{document}

\maketitle
\begin{abstract} We present a new proof, as well as a ${\bf C/Q}$ extension,
 of the
Riemann-Roch-Grothendieck  theorem of Bismut-Lott for flat vector
bundles. The main techniques used are the computations of the
adiabatic limits of $\eta$-invariants associated to the so-called
sub-signature operators. We further show that the Bismut-Lott
analytic torsion form can be derived naturally from the
transgression of the $\eta$-forms appearing  in the adiabatic
limit computations.
\end{abstract}

\renewcommand{\theequation}{\thesection.\arabic{equation}}
\setcounter{equation}{0}

\section{Introduction} \label{s1}

$\quad$ Let $M$ be a compact smooth manifold. For any complex flat
vector bundle $F$ over $M$ with the flat connection $\nabla^F$,
one can define a mod  ${\bf Q}$ version of the Cheeger-Chern-Simons
character $CCS(F,\nabla^F)$ (cf. \cite{CS}) as follows.
Let $k$ be a positive integer such that $kF$ is a topologically
trivial vector bundle. Let $\nabla_0^{kF}$ be a trivial connection on $kF$,
 which can be determined by choosing a global basis of $kF$. Let $k\nabla^F$
 be the connection on $kF$ obtained from the direct sum of $k$
copies of $\nabla^F$. Then we define the mod ${\bf Q}$ version of
the Cheeger-Chern-Simons character as
 \begin{equation}\label{1.1}
CCS(F,\nabla^F)={1\over k}CS(\nabla_0^{kF},k \nabla^{F}),
\end{equation}
 where $CS(\nabla_0^{kF},k \nabla^{F})$ is the Chern-Simons class associated to
 $(kF,k\nabla^F,\nabla_0^{kF})$.  It
 determines a well-defined element in $H^{\rm odd}(M,{\bf C/Q})$
(See Section 2.4 for more details).

 Let $Z\rightarrow M\rightarrow B$ be a fibered manifold with compact base
and fibers. Let $e(TZ)$ be the Euler class of the vertical tangent vector
bundle $TZ$. The flat vector bundle $(F,\nabla^F)$ over $M$ induces
canonically a ${\bf Z}$-graded flat vector bundle
$H^*(Z,F|_Z)=\oplus_{i=0}^{\dim Z}H^i(Z,F|_Z)$ over $B$ (cf. \cite{BL}).
 Let $\nabla^{H^*(Z,F|_Z)}=\oplus_{i=0}^{\dim Z}\nabla^{H^i(Z,F|_Z)}$
 denote the corresponding flat connection induced from $\nabla^F$.

 In \cite{BL}, Bismut and Lott proved a Riemann-Roch-Grothendieck type
formula for the imaginary part of the Cheeger-Chern-Simons character,
which can be stated as an identity in $H^{\rm odd}(B,{\bf R})$,
\begin{equation}\label{1.2}
\int_Z e(TZ){\rm Im}(CCS(F,\nabla^F))=\sum_{i=0}^{\dim Z}(-1)^i{\rm Im}
 (CCS(H^i(Z,F|_Z),\nabla^{H^i(Z,F|_Z)})).
\end{equation}
 They actually proved in \cite{BL} a refinement of (\ref{1.2}) on the differential form level,
 and constructed a real analytic torsion form in the context.

 There is also a topological proof of (\ref{1.2}) given by  Dwyer,
 Weiss and  Williams \cite{DWW}.

 In this paper, we will present a new approach to (\ref{1.2})
 based on  considerations of $\eta$-invariants of
 Atiyah-Patodi-Singer \cite{APS1}. Besides giving a new proof of (\ref{1.2}),
 our method also provides an extension of (\ref{1.2}) to cover
 the real part of the Cheeger-Chern-Simons character.
One of the main results of this paper
 can be stated as the following identity in $H^{\rm odd}(B,{\bf R/Q})$,
 \begin{multline}\label{1.3}
\int_Ze(TZ){\rm Re}(CCS(F,\nabla^F))=\sum_{i=0}^{\dim Z}(-1)^i{\rm Re}
 (CCS(H^i(Z,F|_Z),\nabla^{H^i(Z,F|_Z)})).
\end{multline}
 Putting (\ref{1.2}) and (\ref{1.3}) together, we get the following
formula which can be thought of as a Riemann-Roch-Grothendieck formula
for these Cheeger-Chern-Simons characters.

\begin{thm} \label{t0.1} We have the following identity in
$H^{\rm odd} (B,{\bf C/Q})$,
 \begin{align}\label{1.4}
\int_Ze(TZ)CCS(F,\nabla^F)=\sum_{i=0}^{\dim Z}(-1)^i
 CCS(H^i(Z,F|_Z),\nabla^{H^i(Z,F|_Z)}).
\end{align}
In particular, if ${\bf C}$ denotes the trivial complex line
bundle over $M$, then one has
\begin{align}\label{1.5}
\sum_{i=0}^{\dim Z}(-1)^i
 CCS(H^i(Z,{\bf C}|_Z),\nabla^{H^i(Z,{\bf C}|_Z)})=0 \ \ {
 in}\ \
H^{\rm odd} (B,{\bf C/Q}).
\end{align}
\end{thm}

It turns out that (\ref{1.3}) has been obtained by Bismut in
\cite[Theorem 0.2]{B2}  under the extra condition that $TZ$ is
fiber-wise oriented, while when $\dim Z$ is even, (\ref{1.5}) is
a special case of \cite[Theorem 3.12]{B2}.

 Our proof of (\ref{1.3}), in its full generality, is based on an  extension of
\cite[Theorem 0.2]{Z}, where
 Zhang proved a Riemann-Roch type formula for certain extended versions of
 the Atiyah-Patodi-Singer $\rho$-invariant associated to the sub-signature
 operators constructed also in \cite{Z}. The main method used, as  in \cite{Z}, is the computation
 of the adiabatic limits of the  constructed sub-signature operators, based on
 the techniques developed by Bismut-Cheeger \cite{BC} and Dai \cite{D}, as well as the local
 index computations developed in the papers of Bismut-Lott \cite{BL} and Bismut-Zhang \cite{BZ}.
Moreover, under suitable
  deformations of these sub-signature operators,
the above arguments also lead to a new proof of (\ref{1.2}).
Thus, we  obtain (\ref{1.4}) solely    in the framework of
$\eta$-invariants.

It is particularly interesting that in such a process, the
analytic torsion form of Bismut-Lott \cite{BL} shows up naturally
in a transgression formula of the $\eta$-forms
associated to the deformed operators.
This suggests a natural  relationship between the $\eta$ and
torsion invariants.

We should mention that  the proof in \cite{B2} for  (1.3) relies
also on the computations of adiabatic limits of
$\eta$-invariants. Moreover,
when $\dim Z$ is even, (\ref{1.5})   plays a role in our proof of
(1.3).

From another aspect, in view of the ${\bf R/Z}$-index theory
developed by Lott \cite{L}, one can refine (\ref{1.3})  to an
identity in $K^{-1}_{\bf R/Z}(B)$ if $Z$ is even dimensional and spin$^c$
(cf. Section \ref{s3.19}).

\comment{On the other hand, for any integer $j\geq 0$, Cheeger and
Simons defined in \cite{CS} (see also \cite{B2}) the secondary
character $\widehat{c}_{2j+1}(F,\nabla^F)\in
\widehat{H}^{2j+1}(M,{\bf C/Z})$. Similarly, one has the secondary
characters $\widehat{c}_{2j+1}(H^i(Z,F|_Z),\nabla^{H^i(Z,F|_Z)})$,
$0\leq i\leq \dim Z $, in $\widehat{H}^{ 2j+1} (B,{\bf C/Z})$.  In
view of Theorem \ref{t0.1}, it is nature to formulate

\begin{cj} \label{t0.6} The following identity holds for any integer $j\geq 0$ in
$\widehat{H}^{ 2j+1} (B,{\bf C/Z})$,
 \begin{multline}\label{1.11}
\int_Ze(TZ)\widehat{c}_{2j+1}(F,\nabla^F)=\sum_{i=0}^{\dim
Z}(-1)^i
 \widehat{c}_{2j+1}(H^i(Z,F|_Z),\nabla^{H^i(Z,F|_Z)})\\
 -{\rm rk}(F)\sum_{i=0}^{\dim
Z}(-1)^i
 \widehat{c}_{2j+1}(H^i(Z,{\bf C}|_Z),\nabla^{H^i(Z,{\bf C}|_Z)}).
\end{multline}
\end{cj}}

This paper is organized as follows. In Section \ref{s2}, as in \cite[Section 2]{BL},
 we deal with
  the  finite dimensional situation. In Section \ref{s3}, we develop a proof for both
 (\ref{1.2}) and (\ref{1.3}), and discuss the relations between the $\eta$ and torsion
 forms mentioned above. 

 $\ $

 \noindent {\bf Acknowledgments} The authors would like to thank Professor Jean-Michel Bismut
 for kindly informing us his results in [B3], when we told him our results
 related to (1.3). This work was supported in part by  MOEC and the 973 project
 of MOSTC.
 The authors would like to thank Huitao Feng for
 helpful discussions. Part of the work was done while the second author was visiting the
   Center of Mathematical Sciences
at Zhejiang University in December of 2003, he would like to thank Kefeng Liu and
Hongwei Xu for hospitality.

\section{ \normalsize   $\eta$-invariants and flat cochain complexes}\label{s2}
\setcounter{equation}{0}

$\quad$In this section, we discuss the $\eta$-invariants
associated to a ${\bf Z}$-graded flat cochain complex. The
framework is a combination of those in \cite{BC} and \cite[\S 1,
2]{BL}. We show that the Bismut-Cheeger $\eta$-form is exact in
computing the natural adiabatic limit of $\eta$ invariants
appearing in the context. As a consequence, we deduce an equality
relating the Cheeger-Chern-Simons characters of this cochain
complex and of its cohomology. Moreover, a torsion form is
constructed to transgress the $\eta$-form. This torsion form
turns out to be of the same nature as those constructed by
Bismut-Lott in \cite[\S 2]{BL}.

This section is organized as follows. In Section \ref{s2.1}, we
set up the basic geometric data. In  Section \ref{s2.2}, we
introduce a deformation for the twisted signature operator in the
context. In Section \ref{s2.3}, we  compute the adiabatic limit
of the $\eta$-invariants associated to the deformed twisted
signature operators discussed in  Section \ref{s2.2}. In  Section
\ref{s2.4}, we recall the construction of the mod ${\bf Q}$
Cheeger-Chern-Simons character as well as its relation with
$\eta$-invariants. In  Section \ref{s2.5}, we establish a ${\bf
C/Q}$ formula relating various Cheeger-Chern-Simons characters.
In  Section \ref{s2.6}, we refine the real part of the formula
proved in Section \ref{s2.5} to an identity in $K^{-1}_{\bf
R/Z}$-group. In Section \ref{s2.7}, we construct the torsion form
transgressing the $\eta$-form mentioned above. In  Section
\ref{s2.8}, we discuss in more detail the relationships between
$\eta$ and torsion forms.

 \subsection{\normalsize  Superconnections and flat cochain complexes}
\label{s2.1}

$\quad$Let $(E,v)$ be a ${\bf Z}$-graded cochain complex of
complex vector bundles over a compact smooth manifold $B$,
\begin{align}\label{2.1}
(E,v):\ \ \ 0{\rightarrow} E^0\stackrel{v}{\rightarrow}
E^1\stackrel{v}{\rightarrow}\cdots\stackrel{v}{\rightarrow}
E^n\rightarrow 0.
\end{align}
Let $\nabla^E=\oplus_{i=0}^{n}\nabla^{E^i}$ be a ${\bf Z}$-graded connection on $E$. We call
$(E,v,\nabla^E)$ a flat cochain complex if the following two  identities hold,
\begin{align}\label{2.2}
(\nabla^E)^2=0,\ \ \ \ \ [\nabla^E,v]=0,
\end{align}
where we have adopted the notation of supercommutator in the sense of Quillen \cite{Q}.

Let $h^E=\oplus_{i=0}^nh^{E^i}$ be a ${\bf Z}$-graded Hermitian metric on $E$. Let
$v^*\in C^\infty(B,{\rm Hom}(E^*,E^{*-1}))$ be the adjoint of $v$ respect to $h^E$. Let
$(\nabla^E)^*$ be the adjoint connection of $\nabla^E$ with respect to $h^E$.
Then (cf. \cite[(4.1),(4.2)]{BZ} and \cite[\S 1(g)]{BL})
\begin{align}\label{2.3}
(\nabla^E)^*=\nabla^E+\omega (E,h^E)
\end{align}
with
\begin{align}\label{2.4}
\omega (E,h^E)=(h^E)^{-1}(\nabla^Eh^E).
\end{align}

Let $A'$, $A''$ be the superconnections on $E$ in the sense of Quillen \cite{Q} defined by
\begin{align}\label{2.5}
A'=\nabla^E+v,\ \ \ \ \ \ A''=(\nabla^E)^*+v^*.
\end{align}

Let $N\in {\rm End}(E)$ be the number operator of $E$, i.e., $N$ acts on $E^i$ by
multiplication by $i$. We extend $N$ to an element
of $C^\infty(B,{\rm End}(E))$.

Following \cite[(2.26), (2.30)]{BL}, for any $u>0$, set
\begin{align}\label{2.6}
&C_u'=u^{N/2}A'u^{-N/2}=\nabla^E+\sqrt{u}v,\\
&C_u''=u^{-N/2}A''u^{N/2}=(\nabla^E)^*+\sqrt{u}v^*,\nonumber\\
&C_u={1\over 2}(C_u'+C_u''),\ \ \ \ \ \ \ D_u={1\over 2}(C_u''-C_u'). \nonumber
\end{align}
Let
\begin{align}\label{2.7}
\nabla^{E,e}=\nabla^E+{1\over 2}\omega (E,h^E)
\end{align}
be the Hermitian connection on $(E, h^E)$ (cf. \cite[(1.33)]{BL}
and \cite[(4.3)]{BZ}). Then
\begin{align}\label{2.8}
C_u=\nabla^{E,e}+{\sqrt{u}\over 2}(v+v^*)
\end{align}
is a superconnection on $E$, while
\begin{align}\label{2.9}
D_u={1\over 2}\omega (E,h^E)+{\sqrt{u}\over 2}(v^*-v)
\end{align}
is an element in $C^\infty(B,(\Lambda(T^*B)\widehat{\otimes}{\rm End}(E))^{\rm odd})$.

On the other hand, for any $b\in B$, let $H(E,v)_b=\oplus_{i=0}^nH^i(E,v)_b$ be the cohomology
of the complex $(E,v)_b$. Then as in \cite[\S 2(a)]{BL}, by (\ref{2.2}), there is a ${\bf Z}$-graded
complex vector bundle $H(E,v)$ on $B$ whose fiber over $b\in B$ is $H(E,v)_b$.
Moreover, $H(E,v)$ carries a canonically induced flat connection $\nabla^{H(E,v)}$
(cf. \cite[Proposition 2.5]{BL}).

Also, as in \cite[\S 2(b)]{BL}, it follows from finite dimensional Hodge theory that for
any $b\in B$, there is an isomorphism $H(E,v)_b\simeq \ker((v+v^*)_b)$. Thus, there
is a smooth ${\bf Z}$-graded sub-bundle $\ker(v+v^*)$ of $E$ whose fiber over $b\in B$
is $\ker((v+v^*)_b)$, and
\begin{align}\label{2.10}
H(E,v)\simeq \ker(v+v^*).
\end{align}
As a sub-bundle of $E$, $\ker(v+v^*)$ inherits a Hermitian metric from the Hermitian metric
$h^E$ on $E$. Let $h^{H(E,v)}$ denote the Hermitian metric on $H(E,v)$ obtained via (\ref{2.10}).

Let $p^{\ker(v+v^*)}$ be the orthogonal projection from $E$ onto $\ker(v+v^*)$, it clearly preserves
 the ${\bf Z}$-grading. Then by \cite[Proposition 2.6]{BL}, one knows that
 \begin{align}\label{2.11}
&p^{\ker(v+v^*)}\nabla^Ep^{\ker(v+v^*)}=\nabla^{H(E,v)},\\
& p^{\ker(v+v^*)}\omega(E,h^E)p^{\ker(v+v^*)}=\omega(H(E,v),h^{H(E,v)}),
\nonumber \\
& p^{\ker(v+v^*)}\nabla^{E,e}p^{\ker(v+v^*)}=\nabla^{H(E,v),e}.\nonumber
\end{align}

 \subsection{\normalsize  Twisted signature operators and their deformations}
 \label{s2.2}

$\quad$ We assume in the rest of this section that $p=\dim B$ is
odd and $B$ is oriented.

 Let $g^{TB}$ be a Riemannian metric on $TB$. For $X\in TB$, let $c(X)$, $\widehat{c}(X)$
 be the Clifford actions on $\Lambda(T^*B)$ defined by $c(X)=X^*-i_X$, $\widehat{c}(X)
 =X^*+i_X$, where $X^*\in T^*B$ corresponds to $X$ via $g^{TB}$ (cf. \cite[(3.18)]{BL} and \cite[\S 4(d)]{BZ}).
 Then for any $X$, $Y\in TB$,
 \begin{align}\label{2.12}
& c(X)c(Y)+c(Y)c(X)=-2\langle X,Y\rangle, \\
& \widehat{c}(X)\widehat{c}(Y)+\widehat{c}(Y)\widehat{c}(X)
=2\langle X,Y\rangle,  \nonumber \\
& c(X)\widehat{c}(Y)+\widehat{c}(Y)c(X)=0.\nonumber
\end{align}

 Let $e_1, \cdots, e_p$ be a (local) oriented  orthonormal basis of $TB$.
 Let $N_B$ be the number operator on $\Lambda (T^*B)$. Set
 \begin{align}\label{2.13}
\tau=(\sqrt{-1})^{{p(p+1)\over 2}} (-1)^{N_B+p}
\widehat{c}(e_1)\cdots\widehat{c}(e_p) =(\sqrt{-1})^{{p(p+1)\over
2}}c(e_1)\cdots c(e_p) .
\end{align}
 Then $\tau$ is a well-defined self-adjoint  element such that
 \begin{align}\label{2.14}
\tau^2={\rm Id}|_{\Lambda(T^*B)}.
\end{align}

 Let $\mu$ be a Hermitian vector bundle on $B$ carrying with a
 Hermitian connection $\nabla^\mu$ with the curvature denoted by
$R ^\mu = \nabla^{\mu,2}$.
Let  $\nabla^{TB}$ be the  Levi-Civita connection on $(TB, g^{TB})$ with its curvature $R^{TB}$.
Let $\nabla^{\Lambda(T^*B)}$ be the Hermitian connection on $\Lambda(T^*B)$
 canonically induced from $\nabla^{TB}$.
 Let $\nabla^{\Lambda(T^*B)\otimes\mu\otimes E}$ be the tensor product
 connection on $\Lambda(T^*B)\otimes\mu\otimes E$ given by
 \begin{align}\label{2.15}
\nabla^{\Lambda(T^*B)\otimes\mu\otimes E}=
 \nabla^{\Lambda(T^*B)}\otimes {\rm Id}_{\mu\otimes E}
 +{\rm Id}_{\Lambda(T^*B)}\otimes\nabla^\mu\otimes{\rm Id}_E
 +{\rm Id}_{\Lambda(T^*B)\otimes\mu}\otimes\nabla^E.
\end{align}

 Let the Clifford actions $c$, $\widehat{c}$ extend to actions
 on $\Lambda(T^*B)\otimes\mu\otimes E$ by acting as identity on $\mu\otimes E$.
Let $\varepsilon$ be the induced ${\bf Z}_2$-grading operator on $E$, i.e.,
 $\varepsilon=(-1)^N$ on $E$. We extend $\varepsilon$ to an action on
 $\Lambda(T^*B)\otimes\mu\otimes E$ by acting as identity on $\Lambda(T^*B)\otimes\mu$.

 \begin{defn}\label{t2.1} Let $D_{\rm sig}^{\mu\otimes E}$ be the (twisted)
 signature operator defined by
 \begin{multline}\label{2.16}
D_{\rm sig}^{\mu\otimes E}=\varepsilon\tau\sum_{i=1}^pc(e_i)
 \nabla_{e_i}^{\Lambda^{\rm even}(T^*B)\otimes\mu\otimes E}:
 C^\infty(B,\Lambda^{\rm even}(T^*B)\otimes\mu\otimes E)\\
\rightarrow
 C^\infty(B,\Lambda^{\rm even}(T^*B)\otimes\mu\otimes E).
\end{multline}
\end{defn}
 One verifies  that $D_{\rm sig}^{\mu\otimes E}$ is a formally self-adjoint first order
 elliptic differential operator.

 Let $v$, $v^*$ extend to actions on
 $\Lambda^{\rm even}(T^*B)\otimes\mu\otimes E$ by acting as identity on $\Lambda(T^*B)\otimes\mu$.
 For any $u\geq 0$, set
 \begin{align}\label{2.17}
D_{{\rm sig},u}^{\mu\otimes E}=D_{\rm sig}^{\mu\otimes E}
+{\sqrt{u}\over 2}(v+v^*).
 \end{align}

  \begin{Rem}\label{t2.2}  $D_{{\rm sig},u}^{\mu\otimes E}$ can be thought of as obtained from a
 (signature) quantization of $C_u$. Indeed, if $B$ is spin, then one can
 consider the twisted Dirac operators instead of signature operators.
\end{Rem}

Let $Y_u$ be the skew-adjoint element in ${\rm End}(\Lambda^{\rm even}(T^*B)\otimes\mu\otimes E)$
defined by
\begin{align}\label{2.18}
Y_u={\varepsilon\tau\over 2}\sum_{i=1}^pc(e_i)\omega(E,h^E)(e_i)
+{\sqrt{u}\over 2}(v^*-v).
\end{align}

\begin{defn} \label{t2.3}For any $r\in {\bf R}$ and $u\geq 0$, let
$D_{{\rm sig},u}^{\mu\otimes E}(r)$ be the operator defined by
\begin{multline}\label{2.19}
D_{{\rm sig},u}^{\mu\otimes E}(r)=D_{{\rm sig},u}^{\mu\otimes E}+
\sqrt{-1}r Y_u:C^\infty(B,\Lambda^{\rm even}(T^*B)\otimes\mu\otimes E)\\
\rightarrow
 C^\infty(B,\Lambda^{\rm even}(T^*B)\otimes\mu\otimes E).
\end{multline}
 \end{defn}

 Clearly, $D_{{\rm sig},u}^{\mu\otimes E}(r)$ is still elliptic and formally self-adjoint.

For any $X\in TB$, set
\begin{align}\label{2.20}
 \widetilde{c}(X)=\varepsilon\tau c(X).
\end{align} Then
one verifies that for any $X$, $Y\in TB$,
\begin{align}\label{2.21}
\widetilde{c}(X)\widetilde{c}(Y)+\widetilde{c}(Y)\widetilde{c}(X)=-2
\langle X,Y\rangle.
\end{align}

From (\ref{2.16})-(\ref{2.20}), one deduces that
\begin{multline}\label{2.22}
D_{{\rm sig},u}^{\mu\otimes E}(r)=\sum_{i=1}^k\widetilde{c}(e_i)
 \left(\nabla_{e_i}^{\Lambda^{\rm even}(T^*B)\otimes\mu\otimes E}
 +{\sqrt{-1}r\over 2}
 \omega(E,h^E)(e_i)\right) \\
+{\sqrt{u}\over 2}((1-\sqrt{-1}r)v+(1+\sqrt{-1}r)v^*).
 \end{multline}

 \begin{Rem}\label{t2.4}  One verifies that $\widetilde{c}(X)$,
$X\in TB$, anti-commutes with elements
in ${\rm End}^{\rm odd}(E)$. Thus, we see that we are in a situation closely
related to what  considered in \cite[\S 2(a)]{BC}.
\end{Rem}

\subsection{\normalsize  A computation of adiabatic limits of
 $\eta$-invariants}\label{s2.3}

$\quad$Let $\overline{\eta}(D_{{\rm sig},u}^{\mu\otimes E}(r))$
be the reduced $\eta$-invariant in the sense of
Atiyah-Patodi-Singer \cite{APS1}. More precisely, for
 $s\in \bC,{\rm Re} (s)\geq p$, set
\begin{align}\label{02.23}
\eta (D_{{\rm sig},u}^{\mu\otimes E}(r))(s) = \frac{1}{\Gamma
(\frac{s+1}{2})} \int_0^{+\infty} t^{\frac{s-1}{2}} \tr \left[
D_{{\rm sig},u}^{\mu\otimes E}(r) \exp\left(-tD_{{\rm
sig},u}^{\mu\otimes E}(r) ^2\right)\right]dt.
\end{align}
Then $ \eta (D_{{\rm sig},u}^{\mu\otimes E}(r))(s)$ extends to a
meromorphic function of $s\in\bC$ and is holomorphic at $s=0$. Set
\begin{align}\label{02.24}
\overline{\eta}(D_{{\rm sig},u}^{\mu\otimes E}(r)) = \frac{1}{2}
\left(\eta (D_{{\rm sig},u}^{\mu\otimes E}(r))(0) + \dim \ker
(D_{{\rm sig},u}^{\mu\otimes E}(r)) \right).
\end{align}

By \cite[Theorem 2.7]{BC}, one knows that for any $u\geq 0$,
\begin{align}\label{2.23}
\overline{\eta}(D_{{\rm sig},u}^{\mu\otimes E}(r))\equiv
\overline{\eta}(D_{{\rm sig}}^{\mu\otimes E}(r))\ \ \ {\rm mod}\ \ {\bf Z},
\end{align}
where $D_{{\rm sig}}^{\mu\otimes E}(r)$ is the  notation for
$D_{{\rm sig},u=0}^{\mu\otimes E}(r)$ for brevity.

We fix a square root of $\sqrt{-1}$ and let $\varphi:\Lambda(T^*B)\rightarrow\Lambda(T^*B)$
be the homomorphism defined by $\varphi:\omega\in\Lambda^i(T^*B)\rightarrow
(2\pi \sqrt{-1})^{-i/2}\omega.$ The formulas in what follows will not
depend on the choice of the square root of $\sqrt{-1}$.

Let $\widehat{\eta}_r$ be the $\eta$-form of Bismut-Cheeger \cite[(2.26)]{BC} defined by
\begin{align}\label{2.24}
\widehat{\eta}_r=\left({1\over 2\pi \sqrt{-1}}\right)^{1\over 2}
\varphi \int_0^{+\infty} \tr_s\left[
\left((1-\sqrt{-1}r)v+(1+\sqrt{-1}r)v^*\right)e^{-(C_u+\sqrt{-1}rD_u)^2}
\right]{du\over 4\sqrt{u}},
\end{align}
where $\tr_s$ is the supertrace  on $E$ in the sense of Quillen
\cite{Q} with respect to the $\bZ_2$-grading induced by $(-1)^N$.

\begin{Rem}\label{t2.5} Since $\ker(v+v^*)$ forms a vector bundle over $B$ and
\begin{align}\label{2.25}
\left((1-\sqrt{-1}r)v+(1+\sqrt{-1}r)v^*\right)^2=(1+r^2)(v+v^*)^2,
\end{align}
by \cite[Lemma 2.1]{BC} and \cite[\S 9.1]{BGV}, $\widehat{\eta}_r$ in (\ref{2.24}) is well-defined.
\end{Rem}

Now for any $r\in {\bf R}$, let $D_{{\rm sig}}^{\mu\otimes H(E,v)}(r)$ be the deformed
twisted signature operator defined by replacing $(E,v,\nabla^E, h^E)$
by $(H(E,v),0,\nabla^{H(E,v)},h^{H(E,v)})$. That is,
\begin{align}\label{2.26}
D_{{\rm sig}}^{\mu\otimes
H(E,v)}(r)=\sum_{i=1}^p\widetilde{c}(e_i)
 \left(\nabla_{e_i}^{\Lambda^{\rm even}(T^*B)\otimes\mu\otimes H(E,v)}
 +{\sqrt{-1}r\over 2} \omega(H(E,v),h^{H(E,v)})(e_i)\right).
\end{align}

 \begin{thm}\label{t2.6}  For any $r\in {\bf R}$, the following identity holds,
 \begin{align}\label{2.27}
\overline{\eta}(D_{{\rm sig}}^{\mu\otimes E}(r))\equiv
 \overline{\eta}(D_{{\rm sig}}^{\mu\otimes H(E,v)}(r))\ \ \ {\rm mod}\ \ {\bf Z}.
 \end{align}
\end{thm}

 \begin{proof} By (\ref{2.23}), (\ref{2.27}) is equivalent to
 \begin{align}\label{2.28}
\lim_{u\rightarrow +\infty}\overline{\eta}(D_{{\rm sig},u}^{\mu\otimes E}(r))
 \equiv \overline{\eta}(D_{{\rm sig}}^{\mu\otimes H(E,v)}(r))\ \ \ {\rm mod}\ \ {\bf Z}.
\end{align}

 Now by (\ref{2.25}) and by proceeding as in \cite[Theorem 2.28]{BC},
 one knows  that when $v+v^*$ is invertible, i.e., when
 $H(E,v)=\{0\}$, one has
 \begin{align}\label{2.29}
\lim_{u\rightarrow +\infty}\overline{\eta}(D_{{\rm sig},u}^{\mu\otimes E}(r))
 \equiv \int_B\rL(TB,\nabla^{TB}){\rm ch}(\mu,\nabla^\mu)\widehat{\eta}_r
 ,
\end{align}
 where $\rL(TB,\nabla^{TB})$ is the Hirzebruch characteristic form defined by
 $$\rL(TB,\nabla^{TB})=\varphi
 {\det}^{1/2}\left({R^{TB}\over \tanh\left({R^{TB}/2}\right)} \right),$$
 while ${\rm ch}(\mu,\nabla^\mu)$ is the Chern character form
 defined by
 $$
{\rm ch}(\mu,\nabla^\mu)=\varphi\tr\left[\exp(-R^\mu)\right].$$

While in the general case where $\ker(v+v^*)$ forms a vector
bundle over $B$, one can generalize the arguments in
\cite[Theorem 2.28]{BC} to show that when $\mod \bZ$,
\begin{align}\label{2.30}
\lim_{u\rightarrow +\infty}\overline{\eta}(D_{{\rm sig},u}^{\mu\otimes E}(r))
 \equiv \overline{\eta}(D_{{\rm sig}}^{\mu\otimes H(E,v)}(r))
 +\int_B\rL(TB,\nabla^{TB}){\rm ch}(\mu,\nabla^\mu)\widehat{\eta}_r.
\end{align}

\begin{Rem}\label{t2.7} Indeed, see \cite[Theorem 2.39]{B} for a
very simple proof of (\ref{2.30}).
\end{Rem}

\begin{lemma}\label{t2.8} For any $r\in {\bf R}$,  $\widehat{\eta}_r$
is exact. Moreover,
\begin{align}\label{2.31}
\widehat{\eta}_{r=0}=0.
\end{align}
\end{lemma}
\begin{proof}
From (\ref{2.6}), one verifies directly that for any $r\in {\bf R}$,
\begin{align}\label{2.33}
(C_u+\sqrt{-1}rD_u)^2=(1+r^2)C_u^2=-(1+r^2)D_u^2.
\end{align}

By (\ref{2.8}), (\ref{2.9}),
\begin{align}\label{02.33}
\frac{1}{2\sqrt{u}} (v+v^*) = \frac{-1}{u} [N,D_u],\quad
\frac{1}{2\sqrt{u}} (v^*-v) = \frac{-1}{u} [N,C_u],
\end{align}
from which one gets that for any $r\in{\bf R}$ and $u> 0$,
\begin{multline}\label{2.40}
{1\over 4\sqrt{u}}\tr_s\left[ \left((1-\sqrt{-1}r)v+(1+\sqrt{-1}r)
v^*\right)e^{-(C_u+\sqrt{-1}rD_u)^2} \right]\\
= \frac{-1}{2u} \tr_s \left[([N,D_u]+ \sqrt{-1}r[N,C_u])
e^{-(C_u+\sqrt{-1}rD_u)^2} \right]\\
={\sqrt{-1}r\over
2u}d\tr_s\left[Ne^{-(C_u+\sqrt{-1}rD_u)^2}
\right].
\end{multline}

From (\ref{2.24}) and (\ref{2.40}), one sees that $\widehat{\eta}_r$ is an exact form.
In particular, by setting $r=0$ in (\ref{2.40}) and by (\ref{2.24}), one gets (\ref{2.31}).
\end{proof}

Combining Lemma \ref{t2.8}  with (\ref{2.30}), one gets (\ref{2.28}),
which completes
the proof of Theorem \ref{t2.6}.  \end{proof}

\begin{Rem}\label{t2.9} The transgression formula (\ref{2.40}) suggests that it
is possible to transgress the $\eta$-form $\widehat{\eta}_r$ through torsion like forms
of the same nature as those of Bismut-Lott \cite[Definition 3.22]{BL}. This will be
dealt with in more detail in Section \ref{s2.7}.
\end{Rem}

\subsection{\normalsize   Cheeger-Chern-Simons characters and
$\eta$-invariants}
\label{s2.4}

$\quad$
We first recall the definition of the mod ${\bf Q}$ Cheeger-Chern-Simons character
for flat vector bundles.

Let $W$ be a  complex vector bundle over $B$. Let $\nabla^W_0$,
$\nabla^W_1$ be two connections on $W$. Let $W_t$, $0\leq t\leq
1$, be a smooth path of connections on $W$ connecting $\nabla^W_0$
and $\nabla^W_1$. Let $CS(\nabla^W_0,\nabla^W_1)$ be the
differential form defined by
\begin{align}\label{2.41}
CS(\nabla^W_0,\nabla^W_1)=-\left(1\over 2\pi\sqrt{-1}\right)^{1\over 2}
\varphi\int_0^1\tr\left[{\partial\nabla^W_t\over\partial t}\exp(-(\nabla^W_t)^2)
\right]dt.
\end{align}
Then
\begin{align}\label{2.42}
dCS(\nabla^W_0,\nabla^W_1)={\rm ch}(W,\nabla^W_1)-{\rm ch}(W,\nabla^W_0).
\end{align}
Moreover, up to exact forms, $CS(\nabla^W_0,\nabla^W_1)$ does not depend on the
path of connections
on $W$ connecting $\nabla^W_0$ and $\nabla^W_1$.

\begin{Rem} \label{0t2.9}
If $0\to E ^0 \to E ^1\to E^2\to 0$ is a short exact sequence of
flat vector bundles as in (\ref{2.1}), let $CS(\nabla ^{E^0, e},
\nabla ^{E^1, e}, \nabla ^{E^2, e}) $ be the Chern-Simons class
as in \cite[(10)]{L}, then by Lemma \ref{t2.8} and the
characterization of the   form $\widehat{\eta}$ (cf. \cite[Theorem
2.10]{B1} and
 \cite[Lemma 3.16]{BuMa}), we know that in $\Omega ^{\rm odd} (M)/{\rm Im} (d)$,
\begin{align}\label{02.41}
CS(\nabla ^{E^0, e}, \nabla ^{E^1, e}, \nabla ^{E^2, e}) =
\widehat{\eta}_{r=0} =0.
\end{align}
\end{Rem}

Now let $(F,\nabla^F)$ be a complex flat vector bundle over $B$. Then there is a positive
integer $q$ such that $qF$, the direct sum of $q$ copies of $F$, is topologically
trivial. Let $\nabla_0^{qF}$ be a trivial connection on $qF$ which can be determined
by choosing a global basis of $qF$. Let $q\nabla^F$ be the natural flat connection on $qF$
obtained from the direct sum of $q$ copies of $\nabla^F$.

By (\ref{2.42}), one sees that $CS(\nabla_0^{qF},q\nabla^F)$ is a closed form on $B$.
Moreover, by proceeding as in \cite[Lemma 1]{L}, one shows that
${1\over q}CS(\nabla_0^{qF},q\nabla^F)$ determines a cohomology class in $H^{\rm odd}
(B,{\bf C/Q})$ not depending on the choice of $q$ and $\nabla_0^{qF}$.

\begin{defn}\label{t2.10} We define the mod ${\bf Q}$ Cheeger-Chern-Simons
of $(F,\nabla^F)$ to be
\begin{align}\label{2.43}
CCS(F,\nabla^F)=\left[{1\over q}CS(\nabla_0^{qF},q\nabla^F)\right]\in H^{\rm odd}
(B,{\bf C/Q}).
\end{align}
\end{defn}

By \cite[Proposition 2.9]{CS}, up to ${\rm rk}(F)$,
$CCS(F,\nabla^F)$ is exactly $\widehat{\ch}(F,\nabla^F) \in
\widehat{H}^{\rm odd} (B,{\bf C/Q})$ defined in \cite[(2.19),
Theorem 2.3]{B2}.

Let $h^F$ be a Hermitian metric on $F$. Let $\omega(F,h^F)$ be
given similarly as in (\ref{2.4}), and let $\nabla^{F,e}$ be the
associated Hermitian connection on $F$ given similarly as in
(\ref{2.7}).
 Let $q\nabla^{F,e}$ be
the connection on $qF$ obtained from the direct sum of $q$ copies of
$\nabla^{F,e}$. Then, one verifies directly that
\begin{multline}\label{2.46}
CCS(F,\nabla^F)={1\over q}CS(\nabla_0^{qF},q\nabla^{F,e})+
{1\over q}CS(q\nabla^{F,e},q\nabla^F)\\
={1\over q}CS(\nabla_0^{qF},q\nabla^{F,e})+CS(\nabla^{F,e},\nabla^F).
\end{multline}
Following \cite[(0.2)]{BL},  for any integer $j\geq 0$, let
$c_{2j+1}(F,h^F)$ be the Chern form defined by
\begin{align}\label{2.47}
c_{2j+1}(F,h^F)=(2\pi\sqrt{-1})^{-j}2^{-(2j+1)}\tr\left[\omega^{2j+1}(F,h^F)\right].
\end{align}
Let $c_{2j+1}(F)$ be the associated cohomology class in
$H^{2j+1}(B,{\bf R})$, which does not depend on the choice of
$h^F$. The following identity  has been proved in
\cite[Proposition 1.14]{BL},
\begin{align}\label{2.48}
{1\over \sqrt{-1}}CS(\nabla^{F,e},\nabla^F)
={\rm Im}(CCS(F,\nabla^F))=-{1\over 2\pi}\sum_{j=0}^{+\infty}{2^{2j}j!\over (2j+1)!}
c_{2j+1}(F).
\end{align}
Consequently,
\begin{align}\label{2.49}
{\rm Re}(CCS(F,\nabla^F))={1\over q}CS(\nabla_0^{qF},q\nabla^{F,e})\ \ {\rm in}
\ \ H^{\rm odd}(B,{\bf C/Q}).
\end{align}

We now come to consider the $\eta$-invariants mentioned in the title of this subsection.

Recall that $B$ is compact oriented, carrying with a Riemannian
metric $g^{TB}$ and that $p=\dim B$ is odd. Recall also that
$\mu$ is a Hermitian vector bundle over $B$ carrying with a
Hermitian connection $\nabla^\mu$.

We apply the constructions in the Sections \ref{s2.1}-\ref{s2.3} to the
trivial flat cochain complex $(F,0,\nabla^F)$. Thus, let $D_{\rm
sig}^{\mu\otimes F}$ be the twisted signature operator defined as
in (\ref{2.16}).

Let $h_0^{qF}$ be a Hermitian connection on $qF$ such that $\nabla_0^{qF}$
is a Hermitian connection with respect to $h_0^{qF}$.

It is easy to see that one can construct a smooth pass of Hermitian metrics
connecting $qh^F$ and $h_0^{qF}$, as well as a smooth pass of Hermitian connections
connecting $q\nabla^{F,e}$ and $\nabla_0^{qF}$.

Then by the standard variation formula for reduced $\eta$-invariants (cf. \cite{APS1} and \cite[Theorem 2.10]{BF}),
\begin{align}\label{2.50}
q\overline{\eta}(D_{\rm sig}^{\mu\otimes F})-
q{\rm rk}(F)\overline{\eta}(D_{\rm sig}^{\mu})\equiv\int_B
\rL(TB,\nabla^{TB}){\rm ch}(\mu,\nabla^\mu)CS(\nabla_0^{qF},q\nabla^{F,e})
\ \ {\rm mod}\
{\bf Z}.
\end{align}

From (\ref{2.49}) and (\ref{2.50}), one gets
\begin{align}\label{2.51}
\overline{\eta}(D_{\rm sig}^{\mu\otimes F})-
{\rm rk}(F)\overline{\eta}(D_{\rm sig}^{\mu})\equiv\int_B
\rL(TB){\rm ch}(\mu){\rm Re}(CCS(F,\nabla^{F}))
\ \ {\rm mod}\
{\bf Q}.
\end{align}

For any $r\in{\bf R}$, let $\nabla^{F,e,(r)}$ denote the Hermitian connection
on $F$ defined by
\begin{align}\label{2.52}
\nabla^{F,e,(r)}=\nabla^{F,e}+{\sqrt{-1}r\over 2}\omega(F,h^F).
\end{align}

For any integer $j\geq 0$ and $r\in{\bf R}$, let $a_j(r)\in {\bf
R}$ be defined as
\begin{align}\label{2.53}
a_j(r)=\int_0^1(1+u^2r^2)^jdu.
\end{align}

\begin{lemma}\label{t2.11} The following identity in $H^{\rm
odd}(B,{\bf R})$ holds,
\begin{align}\label{2.54}
CS(\nabla^{F,e},\nabla^{F,e,(r)})=-{r\over
2\pi}\sum_{j=0}^{+\infty} {a_j(r)\over j!}c_{2j+1}(F).
\end{align}
\end{lemma}
\begin{proof} Formula (\ref{2.54}) follows from (\ref{2.41}), (\ref{2.47}) and a direct computation
in considering the smooth pass of connections $(1-u)\nabla^{F,e}+u\nabla^{F,e,(r)}$,
$0\leq u\leq 1$. \end{proof}

\begin{Rem} \label{t2.12} By comparing (\ref{2.48}) and (\ref{2.54}), we see that up to rescaling,
one can recover the imaginary part of the Cheeger-Chern-Simons character
$CCS(F,\nabla^F)$ through (deformed) Hermitian connections.
\end{Rem}

From (\ref{2.22}), (\ref{2.23}), (\ref{2.54}) and the standard variation formula
for reduced $\eta$-invariants, one finds that for any $r\in{\bf R}$,
\begin{multline}\label{2.55}
\overline{\eta}(D_{\rm sig}^{\mu\otimes F}(r))-
\overline{\eta}(D_{\rm sig}^{\mu\otimes F})\equiv\int_B
\rL(TB,\nabla^{TB}){\rm ch}(\mu,\nabla^\mu)
CS(\nabla^{F,e},\nabla^{F,e,(r)})\ \ {\rm mod}\ {\bf Z}\\
=-{r\over 2\pi}\int_B \rL(TB){\rm ch}(\mu)\sum_{j=0}^{+\infty}
{a_j(r)\over j!}c_{2j+1}(F).
\end{multline}

\subsection{\normalsize   Flat cochain complex and the Cheeger-Chern-Simons character}\label{s2.5}

$\quad$ We make the same assumptions and use the same notation as in
Sections \ref{2.1}-\ref{2.3}.
Thus, $(E,v,\nabla^E)$ is a ${\bf Z}$-graded flat cochain complex
over $B$, etc.

Let $CCS(E,\nabla^E)$ denote the ${\bf C/Q}$ Cheeger-Chern-Simons character
defined by
\begin{align}\label{2.56}
CCS(E,\nabla^E)=\sum_{i=0}^n (-1)^iCCS(E^i,\nabla^{E^i})\ \ {\rm in}\ \ H^{\rm odd}
(B,{\bf C/Q}).
\end{align}

The imaginary part of the following result has been proved by Bismut-Lott
\cite[Theorem 2.19]{BL}.

\begin{thm}\label{t2.13} The following identity holds in
$H^{\rm odd} (B,{\bf C/Q})$,
\begin{align}\label{2.57}
CCS(E,\nabla^E)=CCS(H(E,v),\nabla^{H(E,v)}).
\end{align}
 \end{thm}
\begin{proof} By Theorem \ref{t2.6}, we know that for any $r\in{\bf R}$,
\begin{align}\label{2.58}
\overline{\eta}(D_{{\rm sig}}^{\mu\otimes E}(r))-
\overline{\eta}(D_{{\rm sig}}^{\mu\otimes E})\equiv
 \overline{\eta}(D_{{\rm sig}}^{\mu\otimes H(E,v)}(r))
 -\overline{\eta}(D_{{\rm sig}}^{\mu\otimes H(E,v)})\ \ \ {\rm mod}\ \ {\bf Z}.
\end{align}

 From (\ref{2.55}) and (\ref{2.58}), one finds
 \begin{multline}\label{2.59}
{r\over 2\pi}\int_B
\rL(TB){\rm ch}(\mu)\sum_{j=0}^{+\infty}
{a_j(r)\over j!}\sum_{i=1}^n(-1)^ic_{2j+1}(E^i) \\
={r\over 2\pi}\int_B
\rL(TB){\rm ch}(\mu)\sum_{j=0}^{+\infty}
{a_j(r)\over j!}\sum_{i=1}^n(-1)^ic_{2j+1}(H^i(E,v))\ \ {\rm mod}\ {\bf Z}.
\end{multline}
By taking derivative with respect to $r$ at $r=0$, one gets that
\begin{multline}\label{2.60}
\int_B
\rL(TB){\rm ch}(\mu)\sum_{j=0}^{+\infty} {1\over
j!}\sum_{i=1}^n(-1)^ic_{2j+1}(E^i) \\
=\int_B \rL(TB){\rm
ch}(\mu)\sum_{j=0}^{+\infty} {1\over
j!}\sum_{i=1}^n(-1)^ic_{2j+1}(H^i(E,v)).
\end{multline}
Since
(\ref{2.60}) holds for any complex vector bundle $\mu$ over $B$, while
$\rL(TB){\rm ch}(\cdot):K(B)\otimes{\bf Q}\rightarrow
H^{\rm even}(B,{\bf Q})$ is an isomorphism, from (\ref{2.60}) and a
simple degree counting, one deduces that for any integer $j\geq
0$,
\begin{align}\label{2.61}
\sum_{i=1}^n(-1)^ic_{2j+1}(E^i)=\sum_{i=1}^n(-1)^ic_{2j+1}(H^i(E,v))\ \ {\rm in}
\ \ H^{2j+1}(B,{\bf R})  .
\end{align}

From (\ref{2.48}), (\ref{2.56}) and (\ref{2.61}), one gets
\begin{align}\label{2.62}
{\rm Im}(CCS(E,\nabla^E))={\rm Im}(CCS(H(E,v),\nabla^{H(E,v)})),
\end{align}
which was first proved in \cite[Theorem 2.19]{BL} by using a
direct transgression method.

Now by applying (\ref{2.51}) to each $E^i$ as well as each $H^i(E,v)$, $0\leq i\leq n$,
and by Theorem \ref{t2.6}, one finds that
\begin{align}\label{2.63}
\int_B
\rL(TB){\rm ch}(\mu){\rm Re}(CCS(E,\nabla^{E}))\equiv
\int_B \rL(TB){\rm ch}(\mu){\rm
Re}(CCS(H(E,v),\nabla^{H(E,v)})) \ \ {\rm mod}\ {\bf
Q}.
\end{align}
 By using the fact that $\rL(TB){\rm
ch}(\cdot):K(B)\otimes{\bf Q}\rightarrow H^{\rm even}(B,{\bf Q})$
is an isomorphism again, one deduces from (\ref{2.63}) the following
identity in $H^{\rm odd}(B,{\bf R/Q})$,
\begin{align}\label{2.64}
{\rm Re}(CCS(E,\nabla^E))={\rm Re}(CCS(H(E,v),\nabla^{H(E,v)})).
\end{align}

From (\ref{2.62}) and (\ref{2.64}), one gets (\ref{2.57}). \end{proof}

\subsection{\normalsize  A refinement in $K^{-1}_{\bf R/Z}(B)$}
\label{s2.6}

$\quad$ In the discussions in the previous subsections, we have
only assumed that  $B$ is oriented, and this is why we have used
the twisted signature operators. If $B$ is spin$^c$ or even spin,
then we can well use the twisted Dirac operators instead. In
particular, this will enable us to apply the constructions in the
${\bf R/Z}$-index theory  developed by Lott \cite{L} to the
current situation, where the $\eta$-form (at $r=0$) vanishes
tautologically.

In fact, in the language of \cite{L}, one easily sees that
$$(E,\nabla^{E,e},0)
=\sum_{i=1}^n(-1)^i(E^i,\nabla^{E^i,e},0)$$
is an element in $K^{-1}_{\bf R/Z}(B)$.

 \begin{thm}\label{t2.14} The following identity holds in
$K^{-1}_{\bf R/Z}(B)$,
\begin{align}\label{2.65}
(E,\nabla^{E,e},0)
=(H(E,v),\nabla^{H(E,v),e},0).
\end{align}
 \end{thm}
\begin{proof} Clearly, (\ref{2.65}) is a refinement of
the real part of (\ref{2.57}). It is also a direct consequence of
\cite[Def. 6]{L} and Remark \ref{0t2.9}. In fact, let $F^i={\rm
Im} (v_{E^{i-1}}),\  G^i= {\Ker} (v_{E^i})$, then $F^i, G^i$ are
flat vector bundles on $B$ with Hermitian metrics induced by
$h^E$. Now we have short exact sequences of flat vector bundles:
$0\to F^i\to G^i\to H^i(E,v)\to 0$, $0\to G^i\to E^i\to
F^{i+1}\to 0$.
 Then by \cite[Def. 6]{L} and Remark \ref{0t2.9}, in $K^{-1}_{\bf R/Z}(B)$,
\begin{align*}
&(G^i, \nabla^{G^i,e},0)= (F^i, \nabla^{F^i,e},0)+(H^i(E,v), \nabla^{H^i(E,v),e},0), \\
&(E^i, \nabla^{E^i,e},0)=(G^i, \nabla^{G^i,e},0)+ (F^{i+1}, \nabla^{G^{i+1},e},0).
\end{align*}
 Thus we get (\ref{2.65}).
 \end{proof}

 \subsection{\normalsize  Torsion forms and a transgression formula for
 $\widehat{\eta}_r$}
\label{s2.7}

$\quad$
As in \cite[(2.39)]{BL}, we denote
\begin{align}\label{2.68}
d(E)=\sum_{i=0}^n(-1)^ii\, {\rm rk}(E^i),\ \ \
d(H(E,v))=\sum_{i=0}^n(-1)^ii\,  {\rm rk} (H^i(E,v)).
\end{align}

By (\ref{2.33}) and by \cite[Theorem 2.13 and Proposition 2.18]{BL}, one has
that as $u\rightarrow +\infty$,
\begin{align}\label{2.69}
\tr_s\left[Ne^{-(C_u+\sqrt{-1}rD_u)^2}
\right]= d(H(E,v))+O\left({1\over\sqrt{u}}\right),
\end{align}
and that when $u\rightarrow 0^+$,
\begin{align}\label{2.70}
\tr_s\left[Ne^{-(C_u+\sqrt{-1}rD_u)^2}
\right]= d(E)+O(u).
\end{align}

The following definition is closely related to \cite[Definition 2.20]{BL}.
\begin{defn}\label{t2.16}  For any $r\in{\bf R}$, put
\begin{multline}\label{2.71}
I_r=-r  
\varphi
\int_0^{+\infty}
\left(\tr_s\left[Ne^{-(C_u+\sqrt{-1}rD_u)^2}\right]- d(H(E,v)) \right.\\
\left. -(d(E)-d(H(E,v)))e^{-u/4}\right){du\over 2u}.
\end{multline}
\end{defn}

\begin{thm}\label{t2.17}  For any $r\in {\bf R}$, the
 following transgression formula holds,
\begin{align}\label{2.72}
\widehat{\eta}_r=-{1\over 2\pi}dI_r.
\end{align}
\end{thm}

\begin{proof} Formula (\ref{2.72}) follows from
(\ref{2.24}), (\ref{2.40}), (\ref{2.69})-(\ref{2.71}). \end{proof}

Let $T_f(A',h^E)$ be the torsion form constructed in
\cite[Definition 2.20]{BL} associated to the odd holomorphic
function $f(z)$ such that $f'(z)=e^{z^2}$, that is,
\begin{multline}\label{2.73}
T_f(A',h^E)=-
\varphi
\int_0^{+\infty} \left(\tr_s\left[Ne^{D_u^2}\right]- d(H(E,v))
\right.\\
 \left. -(d(E)-d(H(E,v)))e^{-u/4}\right){du\over  2u}.
\end{multline}

\begin{thm}\label{t2.18}   The following identity holds,
\begin{align}\label{2.74}
\left. {\partial I_r\over \partial r}\right|_{r=0}=
T_f(A',h^E).
\end{align}
In particular,
\begin{align}\label{2.75}
\left. {\partial \widehat{\eta}_r\over \partial
r}\right|_{r=0}=-{1\over 2\pi}dT_f(A',h^E).
\end{align}
\end{thm}

\begin{proof} Formula (\ref{2.74}) follows from (\ref{2.71}) and \cite[Definition 2.20]{BL}.
Formula (\ref{2.75}) follows from (\ref{2.72}) and (\ref{2.74}).
\end{proof}

Combining (\ref{2.71}), (\ref{2.72}) with the Bismut-Lott
transgression formula \cite[Theorem 2.22]{BL}, one gets

\begin{cor}\label{t2.19} For any $r\in{\bf R}$, the following identity holds,
\begin{multline}\label{2.77}
\widehat{\eta}_r={r\over 2\pi}\sum_{j=0}^{+\infty}{(1+r^2)^j\over j!}
\sum_{i=1}^n\frac{(-1)^i}{2j+1}c_{2j+1}(H^i(E,v),h^{H^i(E,v)})\\
-{r\over 2\pi}\sum_{j=0}^{+\infty}{(1+r^2)^j\over j!}
\sum_{i=1}^n\frac{(-1)^i}{2j+1}c_{2j+1}(E^i,h^{E^i}).
\end{multline}
In particular,
\begin{multline}\label{2.78}
\left. {\partial\widehat{\eta}_r\over \partial  r}\right|_{r=0}=
{1\over 2\pi}\sum_{j=0}^{+\infty} {1\over
j!}\sum_{i=1}^n\frac{(-1)^i}{2j+1}c_{2j+1}(H^i(E,v),h^{H^i(E,v)})\\
- {1\over 2\pi}\sum_{j=0}^{+\infty} {1\over
j!}\sum_{i=1}^n\frac{(-1)^i}{2j+1}c_{2j+1}(E^i,h^{E^i}).
\end{multline}
\end{cor}

\begin{Rem}\label{t2.20}  In view of Theorems \ref{t2.17} and \ref{t2.18}, a
direct computation of (\ref{2.77}) or (\ref{2.78}) will lead to an alternate
proof of the Bismut-Lott transgression formula \cite[Theorem 2.22]{BL}.
\end{Rem}

\subsection{\normalsize   More on $\eta$ and torsion forms}
\label{s2.8}

$\quad$ On $B\times \bR \times \bR^*_+ $, let $\wi{C} + \sqrt{-1}
r \wi{D}$ be the operator defined by
\begin{align}\label{d7}
(\wi{C} + \sqrt{-1} r \wi{D})_{B\times \{r\}\times\{u\}} = C_u +
\sqrt{-1} r D_u
 + dr \frac{\partial}{\partial r}+ du  \frac{\partial}{\partial u}.
\end{align}
Then $\tr_s \left[ \exp \left(-  (\wi{C} + \sqrt{-1} r
\wi{D})^2\right)\right]$ is closed on $B\times \bR \times \bR^*_+
$, moreover
\begin{multline}\label{d8}
(\wi{C} + \sqrt{-1} r \wi{D})^2 = (C_{u}+\sqrt{-1} r D_{u})^2 \\
+ du \frac{\partial}{\partial u}(C_{u}+\sqrt{-1} r D_{u})
+\sqrt{-1} dr D_{u}.
\end{multline}
By Volterra expansion \cite[\S 2.4]{BGV}, we get
\begin{multline}\label{d9}
\tr_s\left[ \exp \left(-  (\wi{C} + \sqrt{-1} r \wi{D})^2
\right)\right]
= \tr_s\left[\exp \left(-  (C_{u}+\sqrt{-1} r D_{u})^2\right)\right]\\
- du \tr_s\left[\frac{\partial}{\partial u}(C_{u}+\sqrt{-1} r D_{u})
\exp \left(-  (C_{u}+\sqrt{-1} r D_{u})^2\right)\right]\\
-dr \tr_s\left[\sqrt{-1}D_{u}\exp \left(-  (C_{u}+\sqrt{-1} r D_{u})^2\right)\right]\\
+ \int_0^1 ds \, \tr_s\left[du \frac{\partial}{\partial
u}(C_{u}+\sqrt{-1} r D_{u})
\exp \left(- s (C_{u}+\sqrt{-1} r D_{u})^2\right) \right.\\
 \sqrt{-1} dr  D_{u} \exp \left(- (1-s) (C_{u}+\sqrt{-1} r
D_{u})^2\right) \Big].
\end{multline}
Applying the total differentiation $d^{B\times \bR \times
\bR^*_+}$ on $B\times \bR \times \bR^*_+$ to (\ref{d9}), after
comparing the coefficients of $du dr$, and using the fact that
$D_u$ commutes with $\exp \left(- s (C_{u}+\sqrt{-1} r
D_{u})^2\right)$, we get
\begin{multline}\label{d10}
\frac{\partial}{\partial r}\tr_s\left[\frac{\partial}{\partial u}(C_{u}+\sqrt{-1} r D_{u})
\exp \left(-  (C_{u}+\sqrt{-1} r D_{u})^2\right)\right]\\
- \frac{\partial}{\partial u}\tr_s\left[\sqrt{-1}D_{u}\exp \left(-  (C_{u}+\sqrt{-1} r D_{u})^2\right)\right]\\
=\sqrt{-1} d\,  \tr_s\left[\frac{\partial}{\partial
u}(C_{u}+\sqrt{-1} r D_{u})
   D_{u}
\exp \left(-  (C_{u}+\sqrt{-1} r D_{u})^2\right) \right].
\end{multline}
Now by (\ref{02.33}),
\begin{multline}\label{d11}
\tr_s\left[\frac{\partial}{\partial u}(C_{u}+\sqrt{-1} r D_{u})
\sqrt{-1}   D_{u} \exp \left(-  (C_{u}+\sqrt{-1} r D_{u})^2\right)
\right]\\
=-\frac{\sqrt{-1}}{2u} \tr_s\left[( [N_Z, D_u]+ [N_Z,\sqrt{-1} r
C _u])   D_{u} \exp \left(-  (C_{u}+\sqrt{-1} r D_{u})^2\right)
\right]\\
=- \frac{\sqrt{-1}}{2u} \tr_s\left[ 2 N_Z D_{u}^2 \exp \left(-  (C_{u}+\sqrt{-1} r D_{u})^2\right)\right]\\
-\frac{1}{2u} d \tr_s \left[N_Z r D_u \exp \left(-
(C_{u}+\sqrt{-1} r D_{u})^2\right)\right].
\end{multline}

From (\ref{2.40}) and (\ref{d10})-(\ref{d11}), one deduces that
\begin{align}\label{2.83}
{\partial \over \partial r}d\tr_s\left[{rN\over 2u}e^{-(1+r^2)C_u^2}\right]
-{\partial \over \partial u}
\tr_s\left[D_ue^{-(1+r^2)C_u^2}\right]=-{1\over u}
d \tr_s\left[ND_u^2e^{-(1+r^2)C_u^2}\right].
\end{align}
 By taking
$r=0$ in (\ref{2.83}), one gets
\begin{align}\label{2.84}
{\partial \over \partial u}
\tr_s\left[D_ue^{D_u^2}\right]=d\tr_s\left[{N\over
2u}(1+2D_u^2)e^{D_u^2}\right],
\end{align}
 which is exactly \cite[(2.32)]{BL}
 (compare also with \cite[(3.103)-(3.105)]{BL}).

\begin{Rem}\label{t2.21} Formula (\ref{2.84}) plays an essential
role in \cite{BL} in the construction of analytic torsion form. While
it can be proved directly as in \cite{BL}, here we obtain it through
purely considerations of $\eta$-forms. This suggests that there
should be a deep relationship between $\eta$ and torsion
invariants as well as forms.
\end{Rem}

We now come back to (\ref{2.40}). We take derivative with respect to
$r$ in it, and then take $r=0$. What we get is
\begin{align}\label{2.85}
{1\over 4\sqrt{u}}\tr_s\left[ (v^*-v)e^{D_u^2} \right]={1\over
2u}d\tr_s\left[Ne^{D_u^2} \right].
\end{align}
 Together
with (\ref{2.9}), we get
\begin{align}\label{2.86}
\tr_s\left[ {\partial D_u\over\partial u}e^{D_u^2} \right]={1\over
2u}d\tr_s\left[Ne^{D_u^2} \right].
\end{align}

It is interesting to compare (\ref{2.84}) and (\ref{2.86}). In particular, we
can rewrite (\ref{2.86}) as
\begin{align}\label{2.87}
{\partial \over\partial u}\int_0^1\tr_s\left[D_ue^{u^2D_u^2} \right]du={1\over
2u}d\tr_s\left[Ne^{D_u^2} \right].
\end{align}

By (\ref{2.69}), (\ref{2.70}), (\ref{2.75}) and (\ref{2.87}), one can give a direct proof
of (\ref{2.78}). As was pointed out in Remark \ref{t2.20}, this would also
lead to a proof of \cite[Theorem 2.22]{BL}. A comparison like this in
the fibration case would be more interesting.

$\ $

\section{\normalsize Sub-signature operators and a Riemann-Roch
 formula}\label{s3}
\setcounter{equation}{0}

$\quad$In this section, we deal  with the fibration case. We will
give a new proof of the imaginary part of Theorem 1.1, which is a
Riemann-Roch-Grothendieck formula due to Bismut-Lott \cite{BL},
by computing the adiabatic limits of $\eta$ invariants of
deformed sub-signature operators. We will also prove the real
part of Theorem 1.1, in its full generality, by using the same
method. Moreover, we will give a natural derivation of the
Bismut-Lott analytic torsion form \cite{BL} through the
transgression of $\eta$ forms appearing in the adiabatic limit
computations.

This Section is organized as follows. In  Section \ref{s3.1}, we
recall the construction of the Bismut-Lott superconnection
introduced in \cite{BL}. In Section \ref{s3.2}, we define the
sub-signature operator as in \cite{Z}, as well as  a deformation
of this operator. In Section \ref{s3.3}, we state the Lichnerowicz
type formula for the deformed sub-signature operator. In  Section
\ref{s3.4}, we state the main technical result of this section,
Theorem \ref{t3.8}, on the adiabatic limit of the $\eta$
invariants  for the deformed sub-signature operators, which will
be proved in
 Sections  \ref{s3.5} and \ref{s3.6}.
In  Section \ref{s3.9}, we prove Bismut-Lott's formula
(\ref{1.2}) through $\eta$ invariants. In  Section \ref{s3.10}, we
prove (\ref{1.3}) by using
 Theorem \ref{t3.8}. In  Section  \ref{s3.19}, we discuss in details
the relation of our results with Lott's $\bR/\bZ$ index theory \cite{L}.
In  Section  \ref{s3.20}, we will construct the Bismut-Lott
analytic  torsion form through the transgression of $\eta$ forms,
which is one of the main points of view of this paper.

\subsection{\normalsize  The Bismut-Lott Superconnection}\label{s3.1}

$\quad$ Let $\pi:  M\to B$ be a smooth fiber bundle with compact
fiber $Z$ of dimension $n$. We denote by $m=\dim M,\ p=\dim B$.
Let $TZ$ be the vertical tangent bundle of the fiber bundle, and
let $T^* Z$ be its dual bundle.
 Let $F$ be a flat complex vector bundle on $M$ and let $\nabla ^F$
denote its flat connection.

 Let ${TM}=T^HM\oplus TZ$ be a splitting of $TM$.
\comment{
Let $T^H M$ be a sub-bundle of $TM$ such that
\begin{eqnarray}\label{a05}
TM = T^H M \oplus TZ.
\end{eqnarray} }
Let $P^{TZ}, P^{T^HM}$ denote the projection from $TM$ to $TZ, T^HM$.
If $U\in T B$, let $U^H$ be the lift of $U$
 in $T^H M$, so that
$\pi_* U^H = U$.

Let $E= \oplus_{i=0}^n E^i$ be the smooth infinite-dimensional
 {\bf Z}-graded vector bundle over $B$ whose fiber
over $b\in B$ is $C^{\infty}(Z_b, (\Lambda ( T^* Z)\otimes
F)_{|Z_b})$. That is
\begin{align}
C^{\infty}(B, E^i)= C^{\infty}(M, \Lambda^i ( T^* Z)\otimes F).
 \end{align}

\begin{defn} \label{t3.1}
For $s\in C^{\infty}(B, E)$ and $U$ a vector field on $B$, then the Lie
differential $L_{U^H}$ acts on ${C^\infty} (B,E)$. Let $\nabla^E$ be a
{\bf Z}-grading  preserving connection on $E$ defined by
\begin{align}\label{a04}
\nabla^E_U s = L_{U^H} s.
\end{align}
 \end{defn}

If $U_1, U_2$ are vector fields on $B$, put
\begin{align}\label{a03}
T(U_1, U_2)= -P^{TZ} [U_1^H, U_2^H] \in C^\infty (M, TZ).
\end{align}
We denote by $i_T \in \Omega^2 (B, \mbox{Hom} (E^\bullet
,E^{\bullet -1} ))$   the 2-form on $B$ which, to vector fields
$U_1, U_2$ on $ B$,
 assigns the operation of interior multiplication by $T(U_1, U_2)$ on $E$.

Let $d^Z$ be the exterior differentiation along fibers. We
consider $d^Z$ to be an element of $C^{\infty}(B, \mbox{Hom}
(E^\bullet ,E^{\bullet+1} ))$. The exterior differentiation
operator $d^M$, acting on $ \Omega (M,F)=C^\infty(M, \Lambda (T^*
M)\otimes F)$, has degree $1$ and satisfies $(d^M)^2=0$. By
\cite[Proposition 3.4]{BL}, we have
\begin{align}\label{a02}
d^M = d^Z + \nabla^E + i_T.
\end{align}
So $d^M$ is a flat superconnection of total degree $1$ on $E$. We have
\begin{align}\label{a020}
(d^{Z})^2=0, \quad [\nabla^E, d^Z]=0.
\end{align}

Let $g^{TZ}$ be a metric on $TZ$. Let $h^F$ be a Hermitian metric on $F$.
Let $\nabla^{F*}$ be the adjoint of $\nabla^F$ with respect to $h^F$.
Let $\omega(F,h^F)$ and $\nabla^{F,e}$ be the $1$-form on $M$ and
the connection on $F$ defined as in (\ref{2.2}), (\ref{2.7}).

 Let $o(TZ)$ be  the orientation bundle of $TZ$,
a flat  real line bundle on $M$.
Let $dv_Z$ be the Riemannian volume  form on fibers $Z$
associated  to the metric $g^{TZ}$ (Here $dv_Z$ is viewed as a section
of  $\Lambda ^{\dim Z} (T^{*}Z)\otimes o(TZ)$).
Let $\left \langle \  , \  \right \rangle_{\Lambda (T^{*}Z)\otimes F} $
be the metric on $\Lambda (T^{*}Z)\otimes F$ induced by $g^{TZ}, h^F$.
Then $E$ acquires a Hermitian metric $h^E$ such that for
$    \alpha, \alpha' \in C^\infty (B, E)$ and $b \in B$,
\begin{align}\label{a1}
\left \langle \alpha, \alpha' \right \rangle_{h^E} (b)=
\int_{Z_b} \left \langle {\alpha, \alpha' } \right \rangle
_{\Lambda (T^*Z)\otimes F} dv_{Z_b}.
\end{align}

Let $\nabla^{E*}$, $d^{Z*}$, $(d^M)^*$, $(i_T)^*$  be the formal
adjoints of $\nabla^{E}$, $d^Z$, $d^M$, $i_T$ with respect to the
scalar product $\left \langle \, ,  \, \right \rangle_{h^E} $. Set
\begin{align}\label{a2}
&D^Z = d^Z + d^{Z*}, \qquad \quad
\nabla^{E, e}= {1 \over 2} (\nabla ^E + \nabla^{E*}),\\
&\omega(E, h^E) = \nabla^{E*}- \nabla ^E. \nonumber
\end{align}

Let $N_Z$ be the number operator of $E$, i.e. $N_Z$ acts by
multiplication by $k$ on $C^\infty (M, \Lambda^k (T^{*}Z) \otimes
F)$. For $u>0$, set
\begin{align}\label{a3}
&C_u'= u^{N_Z/2} d^M u^{-N_Z/2}, \quad C_u''= u ^{-N_Z/2}(d^M)^* u^{N_Z/2},\\
&C_u = {1 \over 2} (C_u'+C_u''), \quad D_u  = {1 \over 2} (C_u''-C_u'). \nonumber
\end{align}
Then $C_u''$ is the adjoint of $C_u'$ with respect to $h^E$.
 Moreover, $C_u$ is a superconnection on $E$ and $D_u$ is an odd  element of
 $\Omega(B, \mbox{End} (E))$, and
\begin{align}\label{a4}
C_u^2 = - D_u^2,\quad   [C_u, D_u]=0.
\end{align}

Let  $g^{TB}$ be a Riemannian metric on $TB$.
Then $g^{TM}= g^{TZ} \oplus \pi^* g^{TB}$  is a metric on $TM$.
 Let $\nabla^{TM}$, $\nabla^{TB}$
 denote the corresponding
Levi-Civita connections  on $TM, TB$. Put $\nabla^{TZ}= P^{TZ}
\nabla^{TM}$, a connection on $TZ$. As shown in \cite[Theorem
1.9]{B}, $\nabla^{TZ}$ is independent of the choice of $g^{TB}$.
Then ${^0 \nabla} = \nabla^{TZ}\oplus  \pi ^* \nabla^{TB}$ is
also a connection on $TM$. Let $S=\nabla^{TM}- {^0 \nabla}$. By
\cite[Theorem 1.9]{B}, $\left \langle S(\cdot)\cdot,\cdot\right
\rangle_{g^{TM}}$ is a tensor independent of $g^{TB}$. Moreover,
for $U_1, U_2 \in TB$, $X,Y  \in TZ$,
\begin{align}\label{0a4}
&\left\langle S(U_1^H)X, U_2^H \right\rangle_{g^{TM}} =
-\left\langle S(U_1^H)U_2^H, X \right\rangle_{g^{TM}}\\
&\hspace*{25mm}= \left\langle S(X)U_1^H, U_2^H \right\rangle_{g^{TM}}
=\frac{1}{2} \left\langle T(U_1^H, U_2^H), X \right\rangle_{g^{TM}}, \nonumber\\
&\left\langle S(X)Y, U_1^H \right\rangle_{g^{TM}}= - \left\langle
S(X) U_1^H, Y \right\rangle_{g^{TM}} =\frac{1}{2}
(L_{U_1^H}g^{TZ}) (X,Y), \nonumber
\end{align}
and all other terms are zero.

Let $\{f_{\alpha}\}_{\alpha=1}^p$ be an orthonormal basis of $TB$,
set $\{f^{\alpha}\}_{\alpha=1}^p$ the dual basis of $T^*B$. In the
following, it's convenient to identify $f_\alpha$ with
$f_\alpha^H$. Let $\{e_i\}_{i=1}^n$ be an orthonormal basis of
$(TZ, g^{TZ})$. We define a horizontal $1$-form $k$ on $M$ by
\begin{align}\label{a5}
k (f_{\alpha})=- \sum_i \left \langle S(e_i)e_i, f_{\alpha}\right \rangle.
\end{align}
Set
\begin{align}\label{a6}
c(T)= {1 \over 2} \sum_{\alpha,\beta} f^\alpha\wedge f^{\beta}
c \Big (T(f_{\alpha}, f_{\beta}) \Big ),\\
\widehat{c}(T)= {1 \over 2} \sum_{\alpha,\beta} f^\alpha\wedge f^{\beta}
\widehat{c} \Big (T(f_{\alpha}, f_{\beta}) \Big ).\nonumber
\end{align}

Let $\nabla^{\Lambda (T^* Z)}$ be the connection on $\Lambda (T^*
Z)$ induced by $\nabla ^{TZ}$.  Let $\nabla^{TZ\otimes F, e} $ be
the connection on $\Lambda (T^* Z)\otimes F$ induced by
$\nabla^{\Lambda (T^* Z)}$, $\nabla^{F,e}$. Then  by
\cite[(3.36), (3.37), (3.42)]{BL},
\begin{align}\label{a7}
&D^Z= \sum_j  c(e_j)\nabla^{TZ\otimes F, e} _{e_j} - {1 \over 2}\sum_j
\widehat{c}(e_j) \omega(F, h^F)(e_j),\\
&d^{Z*}- d^Z = -\sum_j \widehat{c}(e_j)\nabla^{TZ\otimes F, e} _{e_j}
 + {1 \over 2}\sum_j c(e_j) \omega(F, h^F)(e_j),\nonumber\\
&\nabla^{E, e} =\sum_{\alpha} f^{\alpha}
\left(\nabla^{TZ\otimes F, e}_{f_{\alpha}} + {1 \over 2}
k(f_{\alpha})\right ),\nonumber\\
&\omega(E, h^E) = \sum_{\alpha} f^{\alpha} \Big(\sum_{i,j}\left
\langle S(e_i)e_j, f_{\alpha}\right \rangle  c(e_i)
\widehat{c}(e_j) + \omega(F, h^F)(f_{\alpha})\Big).\nonumber
\end{align}
By \cite[Proposition 3.9]{BL}, we get
\begin{align}\label{a8}
&C_u= {\sqrt{u}  \over 2} D^Z + \nabla ^{E,e} - {1 \over 2\sqrt{u}} c(T),\\
&D_u = {\sqrt{u}  \over 2} (d^{Z*}- d^Z) + {1 \over 2}\omega(E, h^E)
-  {1 \over 2\sqrt{u}} \widehat{c}(T).\nonumber
\end{align}

Let $H^\bullet (Z, F|_Z) = \oplus_{i=0}^{n} H^i (Z, F|_Z) $
be the {\bf Z}-graded vector bundle over $B$ whose fiber over $b\in B$
is the cohomology $H(Z_b,  F_{|Z_b})$ of  the sheaf of locally flat sections
 of $F$ on $Z_b$.
 By \cite[\S 3(f)]{BL}, the flat superconnection $d^M$
 induces a canonical flat connection $\nabla^{H (Z, F|_Z)}$
on $H^\bullet (Z, F|_Z)$ which preserves the {\bf Z}-grading
and which  does not depend on the choice
of $T^H M$.

By Hodge theory, there is an isomorphism
$H^\bullet (Z_b, F_{|Z_b})\simeq \mbox{Ker} (D^{Z_b})$.
Then there is an isomorphism of smooth {\bf Z}-graded vector bundles on $B$
\begin{align}
H^\bullet (Z, F|_Z)\simeq \mbox{Ker} (D^Z).
\end{align}
Clearly $\mbox{Ker} (D^Z)$
 inherits a  metric from the scalar product
$\left \langle \,,\,  \right \rangle_{h^E}  $. Let
$h^{H (Z, F|_Z)}$ be the corresponding metric on $H^\bullet(Z, F|_Z)$.

Let $P$ be the orthogonal projection operator from $E$
on $\mbox{Ker} (D^Z)$ with respect to the Hermitian product (\ref{a1}).
Let $(\nabla^{H (Z, F|_Z)})^*$ be the adjoint of $\nabla^{H (Z, F|_Z)}$
with respect to the Hermitian metric $ h^{H (Z, F|_Z)}$.

The following result is  established in \cite[Proposition 3.14]{BL}.
\begin{prop}\label{t3.2} The following identities hold:
\begin{align}\label{a10}
&\nabla^{H (Z, F|_Z)}= P \nabla^E,\quad
\left (\nabla^{H (Z, F|_Z)}\right  )^*= P \nabla^{E*},\\
&\omega\left (H (Z, F|_Z), h^{H (Z, F|_Z)}\right )= P \omega(E, h^E) P. \nonumber
\end{align}
\end{prop}

\subsection{ \normalsize The sub-signature operator on a fibered manifold}
\label{s3.2}

$\quad$We assume that  $TB$ is oriented.

Let $(\mu, h^\mu)$ be a
Hermitian complex vector bundle over $B$ carrying with a Hermitian
connection $\nabla^\mu$.

Let $N_B, N_M$ be the number operators on $\Lambda (T^*B), \Lambda (T^*M)$,
i.e. they act as multiplication by $k$ on
$\Lambda ^k(T^*B), \Lambda  ^k(T^*M)$ respectively.
Then $N_M=N_B+N_Z$.

Let $\nabla^{\Lambda(T^*M)}$ be the connection on $\Lambda(T^*M)$
canonically induced from $\nabla^{TM}$. Let $\nabla
^{\Lambda(T^*M)\otimes \pi^*\mu\otimes F}$ (resp. $\nabla
^{\Lambda(T^*M)\otimes \pi^*\mu\otimes F,e}$)  be the tensor
product connection on $\Lambda(T^*M)\otimes  \pi^*\mu\otimes F$
induced by $\nabla^{\Lambda(T^*M)}$, $\pi^*\nabla^\mu$ and
$\nabla^{F}$ (resp. $\nabla^{F,e}$).

Let $\{e_a\}_{a=1}^{m}$ be an
 orthonormal basis of $TM$, and its dual basis $\{e ^a\}_{a=1}^{m}$.
 Let
$\{f_\alpha\}_{\alpha=1}^p$  be an oriented orthonormal basis
of $TB$. Set
\begin{align} \label{a17}
&\widehat{\tau} (TB)= (\sqrt{-1})^{p(p+1)\over
2}\widehat{c}(f_1^H)\cdots\widehat{c}(f_p^H),  \\
&\tau (TB)= (\sqrt{-1})^{p(p+1)\over
2} c(f_1^H)\cdots c(f_p^H), \nonumber\\
&\tau = (-1)^{N_Z}\tau (TB). \nonumber
\end{align}
Then the operators $\widehat{\tau} (TB), \tau (TB), \tau$ act naturally on
$\Lambda (T^*M)$, and
\begin{align}\label{a18}
& \widehat{\tau} (TB) ^2= (-1)^p,\quad \tau (TB)^2=\tau ^2 =1,\\
&\tau =  (-1)^p (-1)^{N_M} \widehat{\tau} (TB) =\widehat{\tau}
(TB)(-1)^{N_M}.\nonumber
\end{align}

Let $d^{\nabla^\mu} : \Omega ^a (M, \pi^*\mu\otimes F)\to \Omega
^{a+1} (M, \pi^*\mu\otimes F)$  be the unique extension of
$\nabla ^\mu, \nabla ^F$ which satisfies the Leibniz rule. Let
$d^{\nabla^\mu *} $ be the adjoint of $d^{\nabla^\mu}$ with
respect to the scalar product $\left \langle \ ,\ \right
\rangle_{\Omega (M, \pi^*\mu\otimes F)}$ on $\Omega (M,
\pi^*\mu\otimes F)$ induced by $g^{TM}, h^\mu, h^F$ as in
(\ref{a1}).
 As in \cite[(4.26), (4.27)]{BZ}, we have
\begin{align}\label{a19}
&d^{\nabla^\mu} = \sum_a e^a \wedge \nabla ^{\Lambda(T^*M)\otimes \pi^*\mu\otimes F}_{e_a},\\
&d^{\nabla^\mu *} =- \sum_a i_{e_a} \wedge
\left(\nabla ^{\Lambda(T^*M)\otimes \pi^*\mu\otimes F}_{e_a}
+ \omega (F, h^F)(e_a)\right).\nonumber
\end{align}
For $r\in \bR$,  we introduce the following operators as in
\cite[(1.12)]{Z} and (\ref{2.19}),
\begin{align}\label{a20}
& D^{\pi^*\mu\otimes F}_{\sig}
=\frac{1}{2} \left[ \tau (d^{\nabla^\mu}+d^{\nabla^\mu *})
+ (-1)^{p+1} (d^{\nabla^\mu}+d^{\nabla^\mu *}) \tau \right],\\
&\widehat{D}^{\pi^*\mu\otimes F}_{\sig}
=\frac{1}{2}\left[  \tau (d^{\nabla^\mu *}-d^{\nabla^\mu})
+ (-1)^{p+1} (d^{\nabla^\mu *}-d^{\nabla^\mu}) \tau \right],\nonumber\\
&D^{\pi^*\mu\otimes F}_{\sig}(r)=D^{\pi^*\mu\otimes F}_{\sig}+
\sqrt{-1} r \widehat{D} ^{\pi^*\mu\otimes F}_{\sig}.\nonumber
\end{align}

 Let
$(D^{\pi^*\mu\otimes F}_{\sig})^*$,
$(\widehat{D}^{\pi^*\mu\otimes F}_{\sig})^*$  be the formal
adjoint of $ D^{\pi^*\mu\otimes F}_{\sig}$,
$\widehat{D}^{\pi^*\mu\otimes F}_{\sig}$ with respect to $\left
\langle \ ,\ \right \rangle_{\Omega (M, \pi^*\mu\otimes F)}$. Then
\begin{align}\label{a21}
&\tau D^{\pi^*\mu\otimes F}_{\sig}
= (-1)^{p+1}D^{\pi^*\mu\otimes F}_{\sig}\tau,
 \quad \tau \widehat{D} ^{\pi^*\mu\otimes F}_{\sig}
=(-1)^{p+1}\widehat{D} ^{\pi^*\mu\otimes F}_{\sig}\tau,\\
& (D^{\pi^*\mu\otimes F}_{\sig})^*= (-1)^{p+1}D^{\pi^*\mu\otimes F}_{\sig},
 \quad (\widehat{D}^{\pi^*\mu\otimes F}_{\sig})^*
=(-1)^{p}\widehat{D}^{\pi^*\mu\otimes F}_{\sig}.\nonumber
\end{align}

\begin{Rem}\label{r3.3} If $\mu =\bC$, then $D^{F}_{\sig}$ is different
from the sub-signature operator
 in \cite[(1.12)]{Z} (cf. also \cite{Z1}) by a factor $(\sqrt{-1})^{p(p+1)/2} (-1)^{N_M}$.

 Assume now $M=B, \mu=F=\bC$, then if $p=\dim B$ is odd, $D^{\bC}_{\sig}$ is exactly the
 odd Signature operator in \cite[(2.1)]{APS2}, \cite[(1.38)]{B1}, and $
\widehat{D}^{\bC}_{\sig}=0$; if $p$ is even, then $D^{\bC}_{\sig}= \tau (d+d^*)$ and $\widehat{D}^{\bC}_{\sig}=0$.
\end{Rem}

Following \cite{Z}, we will rewrite
$D^{\pi^*\mu\otimes F}_{\sig}, \widehat{D} ^{\pi^*\mu\otimes F}_{\sig}$
by using the natural connections. Let
$\widetilde{\nabla}^{\Lambda(T^*M)}$ be the Hermitian
connection on $\Lambda(T^*M)$ defined by (cf. \cite[(1.21)]{Z})
\begin{align}\label{a22}
\widetilde{\nabla}^{\Lambda(T^*M)}_X=\nabla^{\Lambda(T^*M)}_X-{1\over
2}\sum_{\alpha=1}^{p}\widehat{c}(P^{TZ}S(X)f_\alpha)\widehat{c}(f_\alpha),\ \ \ X\in TM.
\end{align}

 Let $\widetilde{\nabla}^{e}$ be the tensor product connection on
$\Lambda(T^*M)\otimes  \pi^*\mu\otimes F$ induced
by $\widetilde{\nabla}^{\Lambda(T^*M)}$, $\pi^*\nabla^\mu$
and $\nabla^{F,e}$.  For $r\in \bR$, set
\begin{align}\label{a23}
 &D^{\pi^*\mu\otimes F} = \sum_{a=1}^{m}c(e_a)
\widetilde{\nabla}^{e}_{e_a}
-\frac{1}{2} \sum_{i=1}^{n} \widehat{c}(e_i) \omega(F, h^F)(e_i),\\
 &\widehat{D} ^{\pi^*\mu\otimes F} =- \sum_{i=1}^{n}\widehat{c}(e_i)
\widetilde{\nabla}^{e}_{e_i}
+ \frac{1}{2} \sum_{a=1}^{m} c(e_a) \omega(F, h^F)(e_a)\nonumber\\
&\hspace*{30mm}- \frac{1}{4}\sum_{\alpha,\beta=1}^p \widehat{c}(T(f_\alpha, f_\beta))
\widehat{c}(f_\alpha)\widehat{c}(f_ \beta) ,\nonumber\\
&D^{\pi^*\mu\otimes F}(r) = D^{\pi^*\mu\otimes F}
+ \sqrt{-1} r \widehat{D} ^{\pi^*\mu\otimes F}. \nonumber
\end{align}
The following result extends \cite[Proposition 1.14]{Z}.
\begin{prop} \label{t3.3}
\begin{align}\label{a24}
 &D^{\pi^*\mu\otimes F}_{\sig} = \tau D^{\pi^*\mu\otimes F},\quad
\widehat{D}^{\pi^*\mu\otimes F}_{\sig}=\tau
\widehat{D}  ^{\pi^*\mu\otimes F}.
\end{align}
\end{prop}
\begin{proof} By  (\ref{a19}),
\begin{align}\label{a25}
&d^{\nabla^\mu} + d^{\nabla^\mu *}
= \sum_{a=1}^{m} \left ( c(e_a)\nabla ^{\Lambda(T^*M)\otimes \pi^*\mu\otimes F,e}_{e_a}
-\frac{1}{2}\widehat{c}(e_a)\omega (F, h^F)(e_a) \right ),\\
 &d^{\nabla^\mu *}- d^{\nabla^\mu}
= \sum_{a=1}^{m}
\left ( - \widehat{c}(e_a)\nabla ^{\Lambda(T^*M)\otimes \pi^*\mu\otimes F,e}_{e_a}
+ \frac{1}{2}c(e_a)\omega (F, h^F)(e_a)\right ).
\nonumber
\end{align}
Recall also the following equation was gotten in  \cite[(1.24)]{Z}
by direct computations,
\begin{align}\label{a26}
& \nabla ^{\Lambda(T^*M)}_X \widehat{\tau}(TB)
= - \widehat{\tau}(TB) \sum_{\beta=1}^{p}
\widehat{c}(P^{TZ}S(X)f_\beta)\widehat{c}(f_\beta).
\end{align}
By (\ref{0a4}), (\ref{a17}),
\begin{multline}\label{a27}
\frac{1}{2} (-1)^{p+1}\sum_a \widehat{c}(e_a) \tau \widehat{c}(P^{TZ}S(e_a)f_\beta)\\
= \frac{\tau}{2} \Big( \sum_i
\widehat{c}(e_i)\widehat{c}(P^{TZ}S(e_i)f_\beta)- \sum_\alpha
\widehat{c}(f_\alpha)
\widehat{c}(P^{TZ}S(f_\alpha)f_\beta)\Big)\\
= \frac{\tau}{2} \Big(\sum_i
\widehat{c}(e_i)\widehat{c}(P^{TZ}S(e_i)f_\beta)-\frac{1}{2}
\sum_\alpha \widehat{c}(T(f_\alpha,
f_\beta))\widehat{c}(f_\alpha) \Big).
\end{multline}

Now (\ref{a24}) is a direct consequence of (\ref{a25})-(\ref{a27}).
\end{proof}

From (\ref{a23}), the operator $D^{\pi^*\mu\otimes F}$,
$D^{\pi^*\mu\otimes F}(r)$ are formally self-adjoint
first order elliptic operators,  and
$\widehat{D}^{\pi^*\mu\otimes F}$ is a skew-adjoint
first order differential operator.
 The operator  $D^{\pi^*\mu\otimes F}$ is locally of Dirac type.
By (\ref{a18}), (\ref{a21}) and (\ref{a24}),
\begin{align}\label{a28}
&\tau D^{\pi^*\mu\otimes F} = (-1)^{p+1} D^{\pi^*\mu\otimes F}\tau,\quad
&\tau \widehat{D}^{\pi^*\mu\otimes F} = (-1)^{p+1} \widehat{D}^{\pi^*\mu\otimes F}\tau.
\end{align}

\subsection{\normalsize A Lichnerowicz type formula for
$D^{\pi*\mu\otimes F}_{\sig}(r)$}\label{s3.3}

$\quad$If $B\in {\rm End}(TM)$ is antisymmetric, then the action
of $B$ on $\Lambda (T^* M)$ as a derivation (cf.
\cite[(1.26)]{BGV}) is given by
\begin{eqnarray}\label{b1}
\qquad
\sum_{a,b}\left \langle e_b, B e_a\right \rangle e ^b\wedge i_{e_a}
= {1 \over 4} \sum_{a,b} \left \langle e_b, B e_a \right \rangle
\left(c(e_a)c(e_b)-\widehat{c}(e_i)\widehat{c}(e_b)\right).
\end{eqnarray}
Let $\nabla ^{T^HM}= P^{T^HM}\nabla ^{TM}$ be the connection on
$T^HM$ induced by $\nabla ^{TM}$. Let $R^{TM}, R^{T^HM}, R^{TZ}$
be the curvatures of $\nabla ^{TM}, \nabla ^{T^HM}, \nabla ^{TZ}$
respectively.  Let $K$ be the scalar curvature of $(M, g^{TM})$.
Then
\begin{align*}
\nabla ^{TM} =\nabla ^{T^HM}\oplus \nabla ^{TZ}+ S(\cdot) - P ^{T^HM}S(\cdot)P ^{T^HM}.
\end{align*}
Set
\begin{multline}\label{b3}
 \widehat{R}^e  =-\frac{1}{4}  \sum_{ \alpha ,\beta =1}^{p}
\left\langle R^{T^HM} f_{\alpha}, f_{\beta}\right\rangle
\widehat{c}
(f_{\alpha})\widehat{c}(f_{\beta})\\
-\frac{1}{4}  \sum_{i,j=1}^{n}
\left\langle R^{TZ} e_i,e_j\right\rangle\widehat{c}(e_i)
\widehat{c}(e_j) - \frac{1}{4} \left(\omega\left(F,h^F\right)\right)^2.
\end{multline}
Then by (\ref{a22}), (\ref{b3}), the curvature of $\wi{\nabla} ^e$
is given by
\begin{align}\label{b4}
(\wi{\nabla} ^e)^2  = \frac{1}{4}
\sum_{a,b=1}^m \left\langle R^{TM}e_a,e_b\right\rangle c(e_a)c(e_b) + \widehat{R}^e  + \pi^* R^\mu.
\end{align}
Let $\nabla ^{TM\otimes F,e}\omega(F, h^F)$ be the covariant derivative of $\omega(F, h^F)$. Explicitly
\begin{align}\label{b5}
\nabla ^{TM\otimes F,e}_{e_a}\omega(F, h^F)(e_b) = (\nabla ^{\Lambda (T^*M)\otimes F} _{e_a}\omega(F, h^F) )(e_b)
+\frac{1}{2} (\omega(F, h^F))^2 (e_a,e_b).
\end{align}
Let $\wi{\Delta}^e$ be the Bochner Laplacian
\begin{align}\label{b6}
\wi{\Delta}^e = \sum_{a=1}^m \left((\wi{\nabla}^e_{e_a})^2 -\wi{\nabla}^e_{\nabla^{TM}_{e_a}e_a}\right).
\end{align}

The following result was proved in \cite[Theorem 1.1]{Z} base on
a direct computation.
\begin{prop} \label{t3.4}
\begin{multline}\label{b7}
 D^{\pi^*\mu\otimes F,2} =-\wi{\Delta}^{e}
+\frac{K}{4} + \frac{1}{2}\sum_{a,b=1}^{m}c(e_a)c(e_b)
 (\widehat{R}^e + \pi^* R^\mu) (e_a,e_b) \\
 +{1\over 4}\sum_{i=1}^{n}
\left(\omega\left(F,h^F\right)(e_i)\right)^2
+{1\over 8}\sum_{i,j=1}^{n}\widehat{c}(e_i)\widehat{c}(e_j)
\left(\omega\left(F,h^F\right)\right)^2 (e_i,e_j)\\
-\frac{1}{2} \sum_{a=1}^{m}c(e_a) \Big [\sum_{i=1}^{n}
\widehat{c}(e_i)\nabla ^{TM\otimes F,e}_{e_a}\omega \left(F,h^F\right)(e_i)\\
    + \sum_{\alpha=1}^{p} \widehat{c}(f_\alpha)
 \omega \left(F,h^F\right)(P^{TZ} S(e_a) f_\alpha)\Big ] .
 \end{multline}
\end{prop}
Similarly, for $ \widehat{D}^{\pi^*\mu\otimes F,2}$,
$[D^{\pi^*\mu\otimes F},\widehat{D}^{\pi^*\mu\otimes F}]$, we have
\begin{prop}\label{t3.5}
\begin{multline}\label{b8}
 \widehat{D}^{\pi^*\mu\otimes F,2} =
\sum_{i=1}^n \left( (\wi{\nabla}^{e}_{e_i} )^2
-\wi{\nabla} ^{e}_{\nabla ^{TM}_{e_i}e_i}\right)   + \frac{1}{2}
\sum_{i,j=1}^{n}\widehat{c}(e_i)\widehat{c}(e_j)(\wi{\nabla}^{e})^2
(e_i,e_j)\\
+ \frac{1}{4}\sum_{i=1}^n\widehat{c}(e_i) \Big
[\wi{\nabla}^{e}_{e_i},  \sum_{\alpha,\beta=1}^p
\widehat{c}(T(f_\alpha, f_\beta))
\widehat{c}(f_\alpha) \widehat{c}(f_ \beta)\Big]\\
+\frac{1}{2}\sum_{\alpha,\beta=1}^p  \widehat{c}(f_\alpha) \widehat{c}(f_ \beta) \wi{\nabla}^{e}_{T(f_\alpha, f_\beta)}
- \frac{1}{2}\sum_{i=1}^n\sum_{a=1}^m\widehat{c}(e_i)c(e_a)(\nabla ^{TM\otimes F,e}_{e_i}
\omega \left(F,h^F\right))(e_a)\\
-{1\over 4}\sum_{a=1}^{m}\left(\omega\left(F,h^F\right)(e_a)\right)^2
+{1\over 8}\sum_{a,b=1}^{m}c(e_a)c(e_b)
\left(\omega\left(F,h^F\right)\right)^2 (e_a,e_b)\\
+\frac{1}{16} \Big(\sum_{\alpha,\beta=1}^p \widehat{c}(T(f_\alpha,
f_\beta))\widehat{c}(f_\alpha) \widehat{c}(f_ \beta)\Big)^2,
 \end{multline}
\begin{multline}\nonumber
[D^{\pi^*\mu\otimes F}, \widehat{D}^{\pi^*\mu\otimes F}]=
-\sum_{a=1}^m\sum_{i=1}^n c(e_a)\widehat{c}(e_i) \left(\widehat{R}^e+\pi^* R^\mu
+\frac{1}{4}\omega \left(F,h^F\right)^2\right) (e_a,e_i)\\
-\sum_{\alpha=1}^{p}\omega \left(F,h^F\right)(f_\alpha)\wi{\nabla}^{e}_{f_\alpha}
 + \frac{1}{4}\sum_{\alpha,\beta=1}^p\omega \left(F,h^F\right)( T(f_\alpha, f_\beta))
 \widehat{c}(f_\alpha)\widehat{c}(f_\beta)
\\
+ \frac{1}{4} \sum_{a=1}^mc(e_a) \Big [\wi{\nabla}^{e}_{e_a},
\sum_{\alpha,\beta=1}^p \widehat{c}(T(f_\alpha,
f_\beta))\widehat{c}(f_\alpha) \widehat{c}(f_ \beta)\Big].
 \end{multline}
\end{prop}
\begin{proof}
Note that by (\ref{b5})
\begin{multline}\label{b10}
\nabla ^{TM\otimes F,e}_{e_a} \omega \left(F,h^F\right)(e_b)
- \nabla ^{TM\otimes F,e}_{e_b} \omega \left(F,h^F\right)(e_a)\\
= (\nabla ^{\Lambda (T^*M)\otimes F,e} \omega \left(F,h^F\right)) (e_a,e_b)=0.
 \end{multline}
Thus
\begin{align}\label{b11}
\sum_{a,b=1}^{m}[c(e_a)\wi{\nabla}^e_{e_a}, c(e_b) \omega \left(F,h^F\right)(e_b)]
=- 2\sum_{a=1}^{m}\omega \left(F,h^F\right)(e_a)\wi{\nabla}^e_{e_a}.
\end{align}

By (\ref{0a4}), $\sum_{i=1}^n \widehat{c}(e_i)\widehat{c}(P^{TZ} S(e_i) f_\alpha)
= - k(f_\alpha)$, as in (\ref{b11}), we get
\begin{align}\label{b12}
\sum_{i,j=1}^{m}[\widehat{c}(e_i)\wi{\nabla}^e_{e_i}, \widehat{c}(e_j) \omega \left(F,h^F\right)(e_j)]
=2 \sum_{i=1}^{n}\omega \left(F,h^F\right)(e_i)\wi{\nabla}^e_{e_i}.
\end{align}
We have also
\begin{align}\label{b13}
\left [\widehat{c}(e_i), \widehat{c}(T(f_\alpha, f_\beta))\widehat{c}(f_\alpha) \widehat{c}(f_\beta)\right ]
= 2 \left\langle T(f_\alpha, f_\beta), e_i  \right\rangle\widehat{c}(f_\alpha) \widehat{c}(f_\beta).
\end{align}
From (\ref{a23}), (\ref{b1})-(\ref{b4}), (\ref{b10})-(\ref{b13}),
 and the curvature identity
\cite[Proposition 1.26]{BGV}, we get (\ref{b8}).
\end{proof}

To conclude this subsection, we state the following formula, which
is a consequence of (\ref{a22}) and will be used in a later
occasion.
\begin{multline}\label{b14}
 \left[\wi{\nabla}^e_{e_a}, \widehat{c}(T(f_\alpha, f_\beta))
 \widehat{c}(f_\alpha) \widehat{c}(f_ \beta)\right] = \left[\wi{\nabla}^ {\Lambda (T^* M)}_{e_a}
, \widehat{c}(T(f_\alpha, f_\beta))\widehat{c}(f_\alpha) \widehat{c}(f_ \beta)\right].
\end{multline}

\subsection{\normalsize The $\eta$ invariant for
$D^{\pi*\mu\otimes F}_{\sig}(r)$}\label{s3.4}

$\quad$In this Section, we assume that $B$ is a closed oriented
compact manifold and $p=\dim B$ is odd.

By (\ref{a20}), $D^{\pi*\mu\otimes F}_{\sig}(r)$ preserves the $\bZ_2$-grading on
$\Omega(M,\pi^*\mu\otimes F)$ induced by $(-1)^{N_M}$. We denote
by $D^{\pi*\mu\otimes F}_{\sig,e}(r)$ the restriction of $D^{\pi*\mu\otimes F}_{\sig}(r)$ on
$\Omega ^{\rm even}(M,\pi^*\mu\otimes F)$.
 Let
$\overline{\eta}(D^{\pi*\mu\otimes F}_{\sig,e}(r))$ denote the associated
reduced $\eta$-invariant in the sense of
 \cite{APS1} (cf. (\ref{02.24})). We will omit the notion of $F$ when $F={\bf C}$ is the
 trivial complex line bundle carrying with the trivial metric and
 connection.

\begin{defn}\label{t3.6} Let
$\widetilde{\phi}(M/B,\mu,F,r)\in {\bf R/Z}$ be defined by
\begin{align}\label{b15}
\widetilde{\phi}(M/B,\mu,F,r)
=\overline{\eta}(D^{\pi*\mu\otimes F}_{\sig,e}(r))
-{\rm rk}(F)\overline{\eta}(D^{\pi*\mu}_{\sig,e}(r)) \
{\rm mod}\ {\bf Z}.
\end{align}
\end{defn}
In particular, when $(F,\nabla^F, h^F)$ is unitary,
$\widetilde{\phi}(M/B,\mu,F,0)\in {\bf R/Z}$ is  the
$\rho$-invariant associated to ${D}^{\pi*\mu}_{\sig,e}$ and $F$
in the sense of \cite{APS2}, \cite{APS3}.

\begin{thm}\label{t3.7} i) If $n$ is odd, then
$\overline{\eta}(D^{\pi*\mu\otimes F}_{\sig,e}(r)) \in  {\bf R/Z}$
does not depend on $(g^{TB},g^{TZ},$ $h^\mu,\nabla^\mu)$ and $h^F$.

ii) The number $\widetilde{\phi}(M/B,\mu,F,r)$ does not depend on
$(g^{TB},g^{TZ},h^\mu,\nabla^\mu)$ and $h^F$.
\end{thm}

For any $\var >0$, let $D^{\pi*\mu\otimes F}_{\sig,e,\var}(r)$ be
the operator obtained above by replacing $g^{TB}$ to
$\frac{1}{\var}  g^{TB}$. The following result is the main
technical result of this paper, which generalizes  \cite[Theorem
0.2]{Z}.
\begin{thm}\label{t3.8} We have the following identities  in $\bR/\bZ$,
\begin{align}\label{b16}
\lim_{\var \to 0}
\overline{\eta}(D^{\pi*\mu\otimes F}_{\sig,e,\var}(r))
= \overline{\eta}(D^{\mu\otimes H(Z,F|_Z)}_{\sig,e}(r))
=\sum_{i=0}^{n}(-1)^i\overline{\eta}(D^{\mu\otimes H^i(Z,F|_Z)}_{\sig,e}(r)).
\end{align}
\end{thm}

By applying the previous constructions to the special case with
$M=B$, one  constructs a series
 of smooth invariants
$ \widetilde{\phi}(B,\mu, H^i(Z;F|_Z),r ) $, $0\leq i\leq n$.
They are  the (generalized) $\rho$-invariants associated to
twisted Signature operators on $B$.

\begin{cor}\label{0t3.8} i)  The following identity holds in
${\bf R/Z}$,
\begin{multline}\label{b17}
\widetilde{\phi} (M/B,\mu, F,r)
=\sum_{i=0}^{n}(-1)^i\widetilde{\phi} (B,\mu,H^i(Z;F|_Z),r)\\
-{\rm rk}(F)\sum_{i=0}^{n}(-1)^i\widetilde{\phi} (B,\mu,H^i(Z;{\bf C}|_Z),r).
\end{multline}
ii) If $n=\dim Z$ is odd, then
\begin{align}\label{b18}
\overline{\eta}(D^{\pi*\mu\otimes F}_{\sig,e}(r))
= \overline{\eta}(D^{\mu\otimes H(Z,F|_Z)}_{\sig,e}(r))
 \quad {\rm in}   \, \  \bR/\bZ.
\end{align}
\end{cor}

\subsection{\normalsize  A proof of Theorem \ref{t3.7}}
\label{s3.5}

$\quad$Let $g^{TB}_s, g^{TZ}_s, T^H_sM,  h^F_s, h^\mu_s, \nabla
^\mu_s$ ($s\in \bR$) be a smooth family of the objects as in
Section \ref{s3.1}.

Then for the fibration $\wi{\pi}: \wi{M}=M\times \bR\to  \wi{B}=B\times \bR$,
let $\pi_1: \wi{M}\to M$, $\pi_B:  \wi{B}\to B$ be the natural projections.
 We define $T^H\wi{M}|_{M\times \{s\}}= T^H _sM \oplus \bR$,
$g^ {T\wi{B}}|_{B\times \{s\}} = g^{TB}_s \oplus ds ^2$,
$h^{\pi_{B}^* \mu}|_{B\times \{s\}} =  h^\mu_s$, $h^{\pi_{1}^*
F}|_{B\times \{s\}} =  h^F_s$. Clearly, $\nabla ^{\pi_B^*\mu}=
\nabla ^\mu_s + ds \wedge \frac{\partial}{\partial s}$ is a
Hermitian connection on $(\pi_B^* \mu, h^{\pi_B^* \mu})$.

 We orient $T\wi{B}$ as follow:
if $\{f_\alpha\}_{\alpha=1}^p$ is an oriented orthonormal basis of
$TB$, then the orientation of $T\wi{B}$ is defined by
$f^1\wedge\cdots\wedge f^p \wedge ds$. We denote by $f_{p+1}=
\frac{\partial}{\partial s}$ and
\begin{align}\label{c1}
\wi{\tau} = (-1)^{N_Z} \tau(T\wi{B}).
\end{align}
Now all the construction in Sections \ref{s3.2}, \ref{s3.3} work well for the fibration $\wi{\pi}$.

Let $\Omega_{\pm}(\wi{M}, (\pi\circ\pi_1)^*\mu\otimes \pi_1^*F)$
be the $\pm 1$ eigenspaces of $\wi{\tau}$ in $\Omega (\wi{M},
(\pi\circ\pi_1)^*\mu\otimes \pi_1^*F)$. Then by (\ref{a21}),
(\ref{a28}), $D^{(\pi\circ\pi_1)^*\mu\otimes \pi_1^*F}_{\sig} (r)$
changes the $\bZ_2$-grading induced by $\wi{\tau}$.
\comment{Let
\begin{align*}
D^{(\pi\circ\pi_1)^*\mu\otimes \pi_1^*F}_{\sig, +} (r):
\Omega_{+}(\wi{M}, (\pi\circ\pi_1)^*\mu\otimes \pi_1^*F) \to
\Omega_{-}(\wi{M}, (\pi\circ\pi_1)^*\mu\otimes \pi_1^*F)
\end{align*}
denote the corresponding restriction of
$D^{(\pi\circ\pi_1)^*\mu\otimes \pi_1^*F}_{\sig} (r)$. }

For any $u>0$, let $P_u(x,y)$ be the smooth kernel of $\exp(-u
(D^{(\pi\circ\pi_1)^*\mu\otimes \pi_1^*F} (r))^2)$ with respect
to the Riemannian volume form $dv_{\wi{M}}(y)$. For $x_0\in
\wi{M}$, let $dv_{T_{x_0}\wi{M}}$ be the Riemannian volume form
on $(T_{x_0}\wi{M}, g^{T_{x_0}\wi{M}})$. For $U\in T_{x_0}\wi{M}$,
let $\nabla_U$ be the ordinary derivative in direction $U$. For
$y=(y_1,\cdots,y_{m+1})\in \bR^{m+ 1}$, we identify $y$ as
$\sum_{a=1}^{m+1}y_a e_a$ as a vector in $T_{x_0}\wi{M}$, and set
\begin{multline}\label{c4}
L_{x_0}(r)=
- (1+r^2)\sum_{i=1}^{ n} \left(\nabla_{e_i}
+\frac{1}{4}  \left\langle R^{T\widetilde{M}}_{x_0} y,
e_i\right\rangle   \right)^2\\
-\sum_{\alpha=1}^{p+1} \left(\nabla_{f_\alpha}
+\frac{1}{4}  \left\langle R^{T\widetilde{M}}_{x_0} y,
f_\alpha\right\rangle   \right)^2
-\frac{1}{4}
\sum_{ \alpha ,\beta =1}^{p+1}
\left\langle R^{T^H\wi{M}}_{x_0} f_{\alpha},f_{\beta}\right\rangle
\widehat{c}
(f_{\alpha})\widehat{c}(f_{\beta})\\
-\frac{1}{4} \sum_{i,j=1}^{n}
\left(\left\langle R^{TZ}_{x_0} e_i,e_j\right\rangle
+ r^2 \left\langle R^{T\widetilde{M}}_{x_0} e_i,
e_j\right\rangle  \right)\widehat{c}(e_i)
\widehat{c}(e_j) .
\end{multline}
Let $\exp (-L_{x_0}(r))(y,y')$ ($y,y'\in \bR^{m+1}$)
be the smooth kernel of $\exp (-L_{x_0}(r))$ associated to
$dv_{T_{x_0}\wi{M}}(y')$.

\begin{prop}\label{t3.9} For $x_0\in \wi{M}$,
\begin{multline}\label{c5}
\lim_{u\rightarrow 0}  \tr \left[\wi{\tau} P_u (x_0,x_0)\right]
= (-1)^{\frac{m(m+1)}{2}+p+m} \left(\frac{1}{\pi}\right) ^{\frac{m+1}{2}}  \rk (F)  \\
\int^{\wedge} \widehat{\tau} \left(T\wi{B}\right  )\exp
(-L_{x_0}(r))(0,0) \tr|_\mu \left[\exp\left(- R^{\pi_B^*
\mu}\right)\right],
\end{multline}
here $\int^{\wedge}$ means the coefficient of $e^1\wedge\cdots
\wedge e^{m+1} \, \widehat{c}(e_1)\cdots\widehat{c}(e_{m+1})$ in
$\widehat{\tau} (T\wi{B})$ $\exp (-L_{x_0}(r))(0,0)$$\tr|_\mu
[\exp(- R^{\pi_B^* \mu})]$.
\end{prop}
\begin{proof}  At first, by \cite[Proposition 4.9]{BZ},
among the monomials in terms of
 $c(e_a)$'s and $ \widehat{c}(e_a)$'s,
only $c(e_1)\widehat{c}(e_1)\cdots c(e_{m+1})\widehat{c}(e_{m+1})$
has a nonzero supertrace with the $\bZ_2$-grading on $\Lambda (T^*\wi{M})$ defined by $(-1)^{N_{\wi{M}}}$.
Moreover,
\begin{align}\label{c6}
\tr \left[(-1)^{N_{\wi{M}}}
c(e_1)\widehat{c}(e_1)...c(e_{m+1})\widehat{c}(e_{m+1})\right]
=(-2)^{m+1}.
\end{align}

In view of (\ref{c6}), to compute the local index,
it is convenient to
use the rescaling $\nabla_{e_a}\rightarrow \frac{1}{\sqrt{u}} \nabla_{e_a}$,
$c(e_a)\rightarrow \frac{1}{\sqrt{u}} e^a\wedge -\sqrt{u} i_{e_a}$,
$\widehat{c}(e_a) \rightarrow \widehat{c}(e_a)$ and $y_a\rightarrow \sqrt{u} y_a$.

 We denote by  $L_{1,u}, L_{2,u}, L_{3,u}$ the operators obtained from  $uD^{\pi^*\mu\otimes F,2}$,
$u \widehat{D}^{\pi^*\mu\otimes F,2}$, $u[ D^{\pi^*\mu\otimes F},
\widehat{D}^{\pi^*\mu\otimes F}]$ after the  above rescaling. Then
by  Propositions \ref{t3.4}, \ref{t3.5} and (\ref{b14}),
 as $u\to 0^+$,
\begin{align}\label{c7}
L_{1,u} \to&  - \sum_{a=1}^{m+1} \Big (\nabla_{e_a}
+\frac{1}{4}
\left\langle R^{T\widetilde{M}}_{x_0} y,
e_a\right\rangle \Big)^2 +  \widehat{R}^e_{x_0}
+(\wi{\pi}^* R^{\pi_B^*\mu})_{x_0}.\\
L_{2,u} \to& \sum_{i=1}^{ n} \Big(\nabla_{e_i}
+\frac{1}{4}
\left\langle R^{T\widetilde{M}}_{x_0} y,
e_i\right\rangle  \Big)^2 \nonumber\\
&+\frac{1}{4}
\sum_{i,j=1}^n \left\langle R^{T\widetilde{M}}_{x_0} e_i,e_j
\right\rangle  \widehat{c}(e_i)
\widehat{c}(e_j)
 +\frac{1}{4} \omega \left(\pi_1^*F,h^{\pi_1^*F}\right)^2_{x_0},\nonumber\\
L_{3,u} \to& 0.\nonumber
\end{align}
Thus after rescaling, the operator $u(D^{\pi^*\mu\otimes
F}(r))^2$ has the limit
\begin{align}\label{c8}
L_{x_0}(r)+ (\wi{\pi}^* R^{\pi_B^*\mu})_{x_0}
 - \frac{1+r^2}{4} \omega \left(\pi_1^*F,h^{\pi_1^*F}\right)^2_{x_0}.
\end{align}
By (\ref{c6}), (\ref{c8}) and by proceeding the standard local index technique, we get
\begin{multline}\label{c9}
\lim_{\var\to 0} \tr\left[\wi{\tau} P_u (x_0,x_0)\right]
=(-1)^{\frac{m(m+1)}{2}}
(-2)^{m+1}\left(\frac{1}{4\pi}\right)^{\frac{m+1}{2}} (-1)^{p+1} \\
\int^{\wedge} \widehat{\tau}  \left(T\wi{B}\right)\exp (-L_{x_0}(r))(0,0)\,
\tr|_\mu\left[\exp\left(-R^{\pi_B^* \mu}\right)\right] \\
\tr|_F \left [\exp\left(\frac{1+r^2}{4}\omega \left(\pi_1^*F,h^{\pi_1^*F}\right)^{2}\right)\right ].
\end{multline}
 Now note that (cf. \cite{BL}, (3.77)])
\begin{align}\label{0c9}
\tr|_F \left [\exp\left(\frac{1+r^2}{4}\omega
\left(\pi_1^*F,h^{\pi_1^*F}\right)^{2}\right)\right ]={\rm rk}
(F).
\end{align}

The proof of Proposition \ref{t3.9} is completed.\end{proof}

\begin{proof} ({\it of Theorem  \ref{t3.7}}).
From (\ref{a24}), Proposition \ref{t3.9},
 the Atiyah-Patodi-Singer  index theorem \cite[Theorem 3.10]{APS1},
one get as in \cite[(1.54)]{Z} the following mod $\bZ$ variation
formula of
 $\eta$ invariants,
\begin{multline}\label{c10}
\overline{\eta}\left( D^{\pi^*\mu\otimes F}_{\sig,e,0} (r) \right)
- \overline{\eta}\left( D^{\pi^*\mu\otimes F}_{\sig,e,1} (r)
\right) = (-1)^{\frac{m(m+1)}{2}+p+m}\frac{1}{2}
\left(\frac{1}{\pi}\right)
^{\frac{m+1}{2}}  \rk (F) \\
 \int_{M\times [0,1]} \int^{\wedge} \widehat{\tau}
\left(T\wi{B}\right)\exp
(-L(r))(0,0)\tr|_\mu\left[\exp\left(-R^{\pi_B^*
\mu}\right)\right].
\end{multline}
Now, observe that $L(r)$
has even degree on the Clifford variables $\widehat{c}(e_i)$
(resp. $\widehat{c}(f_\alpha)$), by the parity consideration, we know
 $\int^{\wedge} \widehat{\tau}(T\wi{B})\exp (-L(r))(0,0)$ $
\tr|_\mu[\exp(-R^{\pi_B^* \mu})]$ is zero if $n$ is odd. Thus we
get the first part of Theorem \ref{t3.7}.

On the other hand, when $\dim Z$ is even, from (\ref{c10}) and its
application to the trivial complex line bundle case, we get the
second part of Theorem \ref{t3.7}.

Thus the proof of Theorem \ref{t3.7} is completed.\end{proof}

\subsection{\normalsize  A proof of Theorem \ref{t3.8}}
\label{s3.6}

We will distinguish the objects in Section \ref{s3.2} associated to
 $\frac{1}{\var} g^{TB}$, instead of $g^{TB}$, by adding a subscript $\var$. By (\ref{2.12}),
for $u_1,u_2\in \bR$, ${\bf c}=c$ or $\widehat{c}$,
 \begin{align}\label{c11}
&\var^{N_B/2} {\bf c}_{\var} (u_1 \sqrt{\var} f_\alpha +u_2 e_i))\var^{-N_B/2}
=u_1 {\bf c}(f_\alpha) + u_2 {\bf c}(e_i),\\
& \var^{N_B/2}\tau_\var \var^{-N_B/2}=\tau.\nonumber
\end{align}
Denote by $D^{\pi^*\mu\otimes F}_{s,\var}=
\var^{N_B/2}D^{\pi^*\mu\otimes F}_{\var} \var^{-N_B/2}$,
similarly, we define $\widehat{D} ^{\pi^*\mu\otimes F}_{s,\var}$
and $D^{\pi^*\mu\otimes F}_{s,\var}(r)$.
 Let $\varphi: \Lambda (T^*B)\to \Lambda (T^*B)$ by defined by
$\varphi\omega = (2 \pi \sqrt{-1})^{-\deg \omega /2} \omega$.

Since we have twisted a vector bundle $\mu$ on $B$,  the
superconnection in Section \ref{s3.1} should be modified
accordingly. Let $\nabla ^{E\otimes \mu,e}$ be the connection  on
$E\otimes \mu$ induced by $\nabla ^{E,e}$ and $\nabla^\mu$. Denote
by $C_u^\mu$ the superconnection on $\mu\otimes E$ defined by
replacing $\nabla ^{E,e}$ in (\ref{a8}) by $\nabla ^{E\otimes
\mu,e}$. All other operators in (\ref{a8}) extend naturally on
$E\otimes \mu$. Let $\wi{\nabla}$ be the  connection on
$\Lambda(T^*M)\otimes \pi^*\mu\otimes F$ induced from $\nabla
^{T^HM}\oplus \nabla ^{TZ}$, $\nabla ^\mu$, $\nabla ^{F,e}$.

\begin{thm} \label{t3.10} For any $u>0$, one has
\begin{multline} \label{c12}
{1\over \sqrt{\pi}}\lim_{\var\rightarrow 0} \tr
\left[D^{\pi^*\mu\otimes F}_{\sig,e,\var}(r) \exp\left(-u
D^{\pi^*\mu\otimes F}_{\sig,e,\var}(r)^2\right) \right]\\
 =\int_B
{\rm L}(TB,\nabla ^{TB})
 \ch(\mu, h^\mu)\\
  \left(1\over 2\pi\sqrt{-1}\right)^{1/2}\varphi
 \tr_s\left[2 \sqrt{u}
  \frac{\partial}{\partial u}\left(C_{4u}  + \sqrt{-1}r D_{4u} \right)
 \exp \left(- (1+r^2)C_{4u}^2\right)\right],
\end{multline}
where the $\tr_s$ on $E$ is defined by the $\bZ_2$-grading induced
from $(-1)^{N_Z}$.
\end{thm}
\begin{proof}
Following \cite{BF} and \cite{BC}, let $z$ be an odd Grassmannian
variable which anti-commutes with $c(e_a)$'s and
$\widehat{c}(e_a)$'s. As in \cite[(4.54)]{BC},  if $A$, $B$ are
of trace class in $\End (\Omega^*(M,\pi ^*\mu\otimes F))$, set
\begin{align}\label{c13}
\tr^{z}[A+zB]=\tr[B].
\end{align}
One finds as in \cite{BF} and \cite[(4.55)]{BC} that by
(\ref{a24}), (\ref{a28}), (\ref{c11}),
\begin{multline}\label{c14}
\sqrt{u} \tr \left[D^{\pi^*\mu\otimes F}_{\sig,e,\var}(r)
\exp(-u D^{\pi^*\mu\otimes F}_{\sig,e,\var}(r)^2)\right] \\
= \frac{1}{2}\sqrt{u} \tr \left[\tau_{\var} D^{\pi^*\mu\otimes
F}_{\var}(r)
\exp(-u D^{\pi^*\mu\otimes F}_{\var}(r)^2)\right]\\
=-\frac{1}{2}\tr^{z}\left[\tau_{\var}\exp(-u D^{\pi^*\mu\otimes F}_{\var}(r)^2
 + z\sqrt{u}D^{\pi^*\mu\otimes F}_{\var}(r) )\right]\\
= -\frac{1}{2}\tr^{z}\left[\tau\exp(-u D^{\pi^*\mu\otimes
F}_{s,\var}(r)^2 + z\sqrt{u}D^{\pi^*\mu\otimes F}_{s,\var}(r)
)\right].
\end{multline}

In \cite[Proposition 2.2]{Z}, Zhang formulated a Lichnerowicz type formula
for
$u (D^{\pi^*\mu\otimes F}_{s,\var})^2-z\sqrt{u}D^{\pi^*\mu\otimes F}_{s,\var}$
 which is obtained from (\ref{0a4}) and (\ref{a7}). The corresponding  degenerate term
as $\var \to 0$ is
\begin{align}\label{c15}
-u\var \sum_{\alpha} \Big(\wi{\nabla}_{f_\alpha}
+\frac{\sqrt{\var}}{2} \sum_{i,\beta}\langle S(f_\alpha)e_i,f_\beta \rangle c(e_i)c(f_\beta)
+\frac{zc(f_\alpha)}{2\sqrt{u\var}} \Big) ^2.
\end{align}
On the other hand, from (\ref{0a4}), (\ref{a23}) and Proposition
\ref{t3.5}, it is easy to see that for the operators
$(\widehat{D}^{\pi^*\mu\otimes F}_{s,\var})^2$,
$[D^{\pi^*\mu\otimes F}_{s,\var}, \widehat{D}^{\pi^*\mu\otimes
F}_{s,\var}]$, there is no second order derivative on
$\nabla_{f_\alpha}$ and all the other terms converge as $\var \to
0$. Thus, the only possible singular term in the local index
computation appears in (\ref{c15}).

To cancel this singular term in (\ref{c15}), one can proceed as in
\cite{BC}, \cite{Z}. Here we will  give another argument as in
\cite[\S 7]{BerB}, \cite[\S 7]{Ma}.

We fix $b_0\in B$. For $\delta >0$ small enough, we can identify
the ball $B_0(\delta)\subset T_{b_0}B$ with center $0$ and radius
$\delta$ to the ball in $B$ by using the exponential map. Let
$'\nabla_\var$ be the connection on $\Lambda (\bC(z))
\widehat{\otimes} \Lambda(T^{*}B)$ on $B$ defined by
\begin{align}\label{c16}
{'\nabla}_{\var} ^{\Lambda (\bC(z))
\widehat{\otimes} \Lambda(T^{*}B)}
= \nabla_{\cdot} ^{ \Lambda(T^{*}B)}+ \frac{zc(\cdot)}{2\sqrt{u\var}}.
\end{align}
Then by (\ref{b1}), (\ref{c16}),
\begin{align}\label{c17}
({'\nabla}_{\var})^2  = {1\over 4}\sum_{\alpha,\beta}
\left \langle R^{TB} f_\alpha, f_\beta \right \rangle
\left(c(f_\alpha)c(f_\beta)-\widehat{c}(f_\alpha)\widehat{c}(f_\beta) \right).
\end{align}

Let  $\nabla_{\var}$ be the  connection on $\Lambda
(\bC(z))\widehat{\otimes} \Lambda(T^{*}M) \otimes \pi^*\mu\otimes
F $ $\simeq \Lambda (\bC(z)) \widehat{\otimes}
\Lambda(T^{*}B)\widehat \otimes$ $\Lambda(T^{*}Z) \otimes
\pi^*\mu\otimes F$ induced by $'\nabla_\var$, $\nabla ^{TZ}$,
$\nabla ^{\mu}$ and $\nabla ^{F,e}$.

For $y\in T_{b_0}B$ sufficiently close to $b_0$, we lift
horizontally the path $t\in \bR^*_+ \to ty$ into path $t\in
\bR^*_+ \to x_t\in M$, with $x_t\in Z_{ty}$, ${d x\over dt}\in
T^H M$. For $x_0 \in Z_0$, we identify $(\Lambda (\bC(z))
\widehat{\otimes} \Lambda(T^{*}M)
 \otimes \pi^*\mu\otimes F)_{x_t}$ to  $(\Lambda (\bC(z))
\widehat{\otimes} \Lambda(T^{*}B))_{b_0} \widehat{\otimes}
(\Lambda(T^{*}Z) \otimes \pi^*\mu\otimes F)_{x_0}$
 by parallel  transport along the curve  $t\to x_t \in Z_s$
with respect to the connections  $\nabla_{\var}$.
For $y\in T_{b_0}Y$, set
\begin{align}\label{c18}
{\cal H}=-\sum_{\alpha} \left( \nabla_{f_\alpha}
+\frac{1}{4}\left\langle R^{TB}_{b_0}y, f_\alpha\right\rangle
\right) ^2 - \frac{1}{4} \sum_{\alpha,\beta}\left\langle
R^{TB}_{b_0} f_{\alpha}, f_{\beta}\right\rangle
\widehat{c}(f_{\alpha})\widehat{c}(f_{\beta}) .
\end{align}

Now we do the following Getzler rescaling: $y_\alpha \rightarrow
\sqrt{u\var}y_\alpha$, $\nabla_{f_\alpha}\rightarrow
\frac{1}{\sqrt{u\var}}\nabla_{f_\alpha}$ and
$c(f_\alpha)\rightarrow \frac{1}{\sqrt{u\var}}f^\alpha \wedge
-\sqrt{u\var}i_{f_\alpha}$. By  (\ref{a4}), (\ref{0a4}),
(\ref{a23}), (\ref{b1}), (\ref{c15}),  (\ref{c17}),
\cite[Proposition 2.2]{Z}, and by proceeding similarly as in
\cite[(4.69)]{BC},
 the rescaled operator obtained from
$uD^{\pi^*\mu\otimes F,2}_{s,\var}(r)
-z \sqrt{u} D^{\pi^*\mu\otimes F}_{s,\var}(r)$ converges as
 $\var \to 0$, to
\begin{align}\label{c19}
 {\cal H} + (1+r^2) (C_{4u}^\mu)^2
-z \left(\sqrt{u}D^Z+  \frac{c(T)}{4\sqrt{u}}   +   \sqrt{-1}r
\left(\sqrt{u}(d^{Z*} -d^Z) + \frac{\widehat{c}(T)}{4\sqrt{u}}
\right)\right).
\end{align}
In fact, when $r=0$, this follows from \cite[(2.41)]{Z}; now by
(\ref{0a4}), (\ref{a23}) and Proposition \ref{t3.5}, we find that
 (\ref{c19})  holds in general.

Also from (\ref{a8}),
\begin{align}\label{c20}
(C_u^\mu)^2 = C_u^2 + R^ \mu.
\end{align}
By (\ref{0a4}), (\ref{a8}), (\ref{a18}), (\ref{c6}),
  (\ref{c19}), (\ref{c20}), as in \cite[(2.43)]{Z}, we
 get
\begin{multline}\label{0c21}
\lim_{\var \to 0} \tr \left[D^{\pi^*\mu\otimes F}_{\sig,e,\var}(r)
\exp(-u D^{\pi^*\mu\otimes F}_{\sig,e,\var}(r)^2)\right] \\
=-\frac{1}{2}  (-1)^{\frac{p(p+1)}{2}
+p}\left(\frac{1}{\pi}\right)^p \int_B {\det} ^{1/2} \left(
\frac{R^{TB}/2}{\sinh (R^{TB}/2)}\right)
\tr \left[e ^{-R^\mu}\right]\\
\tr_s \left[2 \sqrt{u}
  \frac{\partial}{\partial u}\left(C_{4u}  + \sqrt{-1}r D_{4u} \right)
 \exp \left(- (1+r^2)C_{4u}^2\right)\right]\\
 \int^{\wedge_B} \widehat{\tau} (TB)  \exp\left(\frac{1}{4}\sum_{\alpha,\beta} \left\langle R^{TB}
 f_\alpha, f_\beta
\right\rangle  \widehat{c}(f_\alpha)\widehat{c}(f_\beta) \right),
\end{multline}
where $\int^{\wedge_B}$ means the coefficient of $f^1\cdots
f^p\wedge\widehat{c}(f_1)\cdots \widehat{c}(f_p)$ in the last term
of (\ref{0c21}). From (\ref{a17}), (\ref{0c21}), as in
\cite[(2.44)]{Z}, the last term of (\ref{0c21}) is
$(\sqrt{-1})^{\frac{p(p+1)}{2}}$ $\det ^{1/2} (\cosh (R^{TB}/2))$.
Thus we get (\ref{c12}).
\end{proof}

The following Lemma tells us that the right hand side of
(\ref{c12}) is zero.
\begin{lemma}\label{t3.11}
\begin{multline}\label{c21}
\tr_s\left[
\left(\frac{\partial}{\partial u} \left(C_{u}+ \sqrt{-1}r  D_{u}\right) \right)
\exp \left(-  (C_{u}+\sqrt{-1} r D_{u})^2\right)\right]\\
= \frac{\sqrt{-1}r}{2u} d\, \tr_s\left[ N_Z \exp \left(-  (C_{u}+\sqrt{-1} r D_{u})^2\right)\right].
\end{multline}
\end{lemma}
\begin{proof} By (\ref{a4}),
\begin{align}\label{c22}
(C_{u}+\sqrt{-1} r D_{u})^2 = (1+r^2)C_{u}^2.
\end{align}
 By (\ref{a02}), (\ref{a3}), we have (cf. also \cite[p19]{Ma})
\begin{align}\label{c23}
2u \frac{\partial}{\partial u} C_u = - [N_Z, D_u], \quad
2u \frac{\partial}{\partial u} D_u = - [N_Z, C_u].
\end{align}
Thus by  (\ref{c22}), (\ref{c23})
\begin{multline}\label{c24}
\tr_s\left[\frac{\partial}{\partial u} \left(C_{u}+ \sqrt{-1}r  D_{u}\right)
\exp \left(-  (C_{u}+\sqrt{-1} r D_{u})^2\right)\right]\\
= \frac{-1}{2u} \tr_s\left[\left( [N_Z, D_u]+ \sqrt{-1}r[N_Z, C_u]\right)
\exp \left(-  (C_{u}+\sqrt{-1} r D_{u})^2\right)\right]\\
=\frac{\sqrt{-1}r}{2u}  \tr_s\left[ \left[C_u, N_Z \exp \left(-  (C_{u}+\sqrt{-1} r D_{u})^2\right)\right]\right]\\
= \frac{\sqrt{-1}r}{2u} d\, \tr_s\left[ N_Z \exp \left(-  (C_{u}+\sqrt{-1} r D_{u})^2\right)\right].
\end{multline}
\end{proof}

\begin{proof} ({\it of Theorem \ref{t3.8}}). From Lemma \ref{t3.11}, we know that
the right side of (\ref{c12}) is zero. By mimicing the argument
in \cite{D}, we get Theorem \ref{t3.8} (compare also with
\cite[p. 298]{L} for the precise counting of the $
\mod \bZ$ term).
\end{proof}

\subsection{\normalsize $\eta$ invariant and Bismut-Lott theorem}
\label{s3.9}

We prove the Bismut-Lott formula (1.2) in this subsection.

 By the
standard variation formula for reduced eta invariants, we find
that for any $r \in\bR$,
\begin{multline}\label{f1}
\frac{\partial}{\partial r} \overline{\eta}(D^{\pi*\mu\otimes F}_{\sig,e}(r))\\
= \lim_{u\to 0} \sqrt{u} \tr \left[ \Big(\frac{\partial}{\partial r} D^{\pi*\mu\otimes F}_{\sig,e}(r)\Big)
 \exp\left(-u(D^{\pi*\mu\otimes F}_{\sig,e}(r))^2\right)\right]\\
=\lim_{u\to 0}  \tr \left[ \sqrt{-1}\sqrt{u} \widehat{D}^{\pi*\mu\otimes F}_{\sig,e}
 \exp\left(-u(D^{\pi*\mu\otimes F}_{\sig,e}(r))^2\right)\right]\\
=\lim_{u\to 0} \frac{\sqrt{-1}}{2}  \tr \left[\tau \sqrt{u} \widehat{D}^{\pi*\mu\otimes F}
\exp\left(-u(D^{\pi*\mu\otimes F}(r))^2\right)\right].
\end{multline}
Now we do the rescaling as in Section \ref{s3.5}, then by
(\ref{a23}), we know the rescaled operator of $\sqrt{u}
\widehat{D}^{\pi*\mu\otimes F}$ is
\begin{align}
\sum_{i=1}^n\widehat{c}(e_i)\Big(\nabla_{e_i} +\frac{1}{4}
\left\langle R^{T\widetilde{M}}_{x_0} y, e_i\right\rangle \Big)
+\frac{1}{2} \omega (F, h^F)_{x_0}.
\end{align}
By (\ref{c4}),
(\ref{c7}), the rescaled operator of $u(D^{\pi*\mu\otimes
F}(0))^2$ is $L_{x_0}(0)$, and
\begin{align}
\sum_{i=1}^n \widehat{c}(e_i)\Big(\nabla_{e_i} +\frac{1}{4}
\left\langle R^{T\widetilde{M}}_{x_0} y, e_i\right\rangle
\Big)\exp (-L_{x_0}(0))(0,0)
\end{align}
is zero.  Thus as in Section
\ref{s3.5}, by (\ref{c6}), we get
\begin{multline}\label{f2}
\lim_{u\to 0} \frac{\sqrt{-1}}{2}  \tr \left[\tau\sqrt{u} \widehat{D}^{\pi*\mu\otimes F}
\exp\left(-u(D^{\pi*\mu\otimes F}(0))^2\right)\right]\\
=\frac{1}{2} (-1)^{\frac{m(m+1)}{2}+p+m}
\left(\frac{1}{\pi}\right)^{\frac{m}{2}} \sqrt{-1} \int_M dv_M
\int^{\wedge} \widehat{\tau} (TB)\exp (-L_{x_0}(0))(0,0) \\
\tr|_\mu\left[\exp\left(-R^{\pi ^* \mu}\right)\right] \tr|_F \left
[ \frac{1}{2} \omega (F, h^F) \exp\left(-\frac{1}{4}\omega
\left(F,h^{F}\right)^{2}\right)\right ].
\end{multline}
By a simple algebraic result (cf. \cite[Prop. 3.13]{BGV}), the
Berezin integral of $\int ^{\wedge}$ in (\ref{f2})  is (cf.
\cite[(1.49)-(1.51)]{Z}) the coefficient of $\widehat{e ^1}\cdots
\widehat{e ^n}$ in
\begin{multline*}
(\sqrt{-1}) ^{\frac{p(p+1)}{2}} \tr|_\mu\left[\exp\left(-R^{\pi ^*
\mu}\right)\right] \tr|_F \left [ \frac{1}{2} \omega (F, h^F)
\exp\left(-\frac{1}{4}\omega
\left(F,h^{F}\right)^{2}\right)\right ]\\
{\det}^{1/2} \left(\frac{R^{TM}/2}{\sinh{R^{TM}/2}}\right)
{\det}^{1/2}\left( \cosh\left({R^{T^H M}\over 2}\right) \right)
{\det}^{1/2}
\left(\frac{\sinh(R^{TZ}/2)}{R^{TZ}/2}\right)\\
\exp\Big(\frac{1}{4} \sum_{i,j} \left\langle R^{TZ}
e_i,e_j\right\rangle \widehat{e ^i} \wedge \widehat{e ^j} \Big).
\end{multline*}
Thus at $r=0$, we get
\begin{multline}\label{f3}
\frac{\partial}{\partial r} \overline{\eta}(D^{\pi*\mu\otimes
F}_{\sig,e}(r)) _{r=0} =  -\frac{1}{2}\int_B {\rm L}(TB) \ch
(\mu)\int_Z e(TZ)
 \sum_{j=0}^\infty \frac{1}{j!} c_{2j+1}(F).
\end{multline}

By applying (\ref{f3}) to the case of $M=B$ and by (\ref{2.55})
and  Theorems \ref{t3.7}, \ref{t3.8},  one finds,
\begin{multline}\label{f-1}\int_B
{\rm L}(TB) \ch (\mu)\int_Z e(TZ) \sum_{j=0}^\infty \frac{1}{j!} \int_Z e(TZ)  c_{2j+1}(F)\\
=\int_B{\rm L}(TB) \ch (\mu)\sum_{j=0}^\infty \frac{1}{j!}
c_{2j+1}(H(Z, F|_Z)).
\end{multline}

As ${\rm L}(TB)\ch(\cdot): K(B)\otimes \bR\to H^{\rm even}
(B,\bR)$ is an isomorphism, we get from (\ref{f-1}) that,
 \begin{align}\label{f4}
\sum_{j=0}^\infty \frac{1}{j!}  \int_Z e(TZ)  c_{2j+1}(F) =
\sum_{j=0}^\infty \frac{1}{j!} c_{2j+1}(H(Z, F|_Z))  \quad  {\rm
in} \,  H^{\rm odd} (B,\bR) ,
\end{align}
which is equivalent to the Bismut-Lott formula (1.2) through a
simple degree counting, in the case where $B$ is orientable and
of odd dimension.

\subsection{\normalsize  A proof of Theorem \ref{t0.1}  }
\label{s3.10}

$\quad$First of all, if $B$ is not orientable, then there is a
double covering $\sigma: B'\to B$ such that $B'$ is orientable. We
pull-back the fibration $\pi:M\to B$ to get a fibration $\pi' :
M'\to B'$. Thus we only need to prove Theorem \ref{t0.1} when $B$
is orientable. And from now on, we assume  $B$ is orientable.

Next, if $B$ is of even dimension, then we can apply the analysis
before to the product fibration $Z\rightarrow M\times
S^1\rightarrow B\times S^1$ to get the result. Thus from now on,
we can also assume that $\dim B$ is odd.

Combining with what was done in the last subsection, we get a new
proof of Bismut-Lott formula (1.2).

It remains to prove (\ref{1.3}).

As in (\ref{2.51}), we have, when  $\mod \bQ$,
\begin{align}\label{f5}
\widetilde{\phi} (M/B,\mu, F,0)  = \int_B {\rm L}(TB) \ch(\mu)
\int_Z e(TZ) {\rm Re} (CCS (F,\nabla ^F)).
\end{align}
Thus from Theorem \ref{t3.8}, (\ref{2.51}), (\ref{f5}), and the
argument as in the proof of (\ref{f4}), we get in $H^{\rm
odd}(B,{\bf R/Q})$,
 \begin{multline}\label{f6}
\int_Ze(TZ){\rm Re}(CCS(F,\nabla^F))=\sum_{i=0}^{n}(-1)^i{\rm Re}
 (CCS(H^i(Z,F|_Z),\nabla^{H^i(Z,F|_Z)}))\\
-{\rm rk}(F)\sum_{i=0}^{n}(-1)^i
 {\rm Re}(CCS(H^i(Z,{\bf C}|_Z),\nabla^{H^i(Z,{\bf C}|_Z)})).
\end{multline}

Now by Poincar\'e duality, one has for any nonnegative integer
$i$,
\begin{align}\label{z1}
H^i(Z,(F\otimes o(TZ))|_Z)=(H^{n-i}(Z,F^*|_Z))^*.
\end{align}
Thus, if $F=\overline{F^*}$, one has
\begin{multline}\label{z2}
\sum_{i=0}^{n}(-1)^i{\rm Re}
 (CCS(H^i(Z,({F}\otimes o(TZ))|_Z),\nabla^{H^i(Z,({F}\otimes o(TZ))|_Z)}))\\
 =
 \sum_{i=0}^{n}(-1)^{n-i}{\rm Re}
 (CCS((H^i(Z,{F^*}|_Z))^*,\nabla^{(H^i(Z,{F^*}|_Z))^*}))\\
 =\sum_{i=0}^{n}(-1)^{n-i}{\rm Re}
 (CCS((\overline{H^i(Z,{F}|_Z)})^*,\nabla^{(\overline{H^i(Z,{F}|_Z)})^*}))\\
 =(-1)^{n}\sum_{i=0}^{n}(-1)^{i}{\rm Re}
 (CCS(H^i(Z,F|_Z),\nabla^{H^i(Z,F|_Z)})).
\end{multline}

If $n=\dim Z$ is odd, then by setting $F={\bf C}\otimes o(TZ)$ in
(\ref{f6}) and by (\ref{z2}), one gets
\begin{align}\label{z3}
2\sum_{i=0}^{n}(-1)^i {\rm Re} (CCS(H^i(Z,{\bf
C}|_Z),\nabla^{H^i(Z,{\bf C}|_Z)}))=0\ \ {\rm in}\ \ H^{\rm
odd}(B,{\bf R/Q}),
\end{align}
from which (\ref{1.5}) follows.

On the other hand, if $n$ is even, then (\ref{1.5}) follows from
the second part of \cite[Theorem 3.12]{B2}.

Thus, (\ref{1.5}) holds in its full generality.

From (\ref{1.5}) and (\ref{f6}), one gets (\ref{1.3}).

The proof of Theorem \ref{t0.1} is completed.

\subsection{ \normalsize A refinement in $K^{-1}_{\bR/\bZ} (B)$}
\label{s3.19}

$\quad$
Recall that
 $\pi : M\to B$ is a fibration of compact smooth manifolds with compact
fiber $Z$, and   $g ^{TZ}$ is a metric on $TZ$.
 Let  $(F,\nabla ^F)$ be a complex flat vector bundle and
$h ^F$ is a Hermitian metric on $F$. Recall also that
 $\nabla ^{F,e}$ is the Hermitian connection on $F$ induced by
$\nabla^F$ and $h^F$ as in (\ref{2.7}).

 Suppose that $Z$ is even dimensional and spin$^c$.
Let $S(TZ)= S^+(TZ)\oplus S^-(TZ)$ be the spinor bundle of $TZ$.
In \cite[\S 4]{L}, Lott defined a topological index $\ind_{top}$
and an analytic index $ \ind_{an}$, mapping from $
K^{-1}_{\bR/\bZ} (M)$ to $ K^{-1}_{\bR/\bZ} (B)$. Especially, for
any $\mathcal{E}\in K^{-1}_{\bR/\bZ} (M)$ \cite[(37)]{L},
\begin{align}\label{h1}
\ch_{\bR/\bQ} (\ind_{top}(\mathcal{E}))= \int_Z \widehat{A}(TZ)
 e ^{c_1(L_Z)/2}\ch_{\bR/\bQ}(\mathcal{E}),
\end{align}
where $\widehat{A}(TZ)= \varphi
 {\det}^{1/2}\left({R^{TZ}\over \sinh\left({R^{TZ}/2}\right)} \right)$
is the $A$-hat genus of $TZ$ and $c_1(L_Z)$ is the first Chern
class of the complex line bundle $L_Z$ which defines the spin$^c$
structure of $TZ$. The main result of Lott \cite[Corollaries  1
and  3]{L} is that for any $\mathcal{E}\in K^{-1}_{\bR/\bZ} (M)$,
\begin{align}\label{h2}
\ind_{top}(\mathcal{E})= \ind_{an}(\mathcal{E}).
\end{align}

We denote  by $\bC$ the trivial
complex line bundle carrying with the trivial metric and connection.
Then $\mathcal{F}= [(F, h^F, \nabla ^{F,e}, 0)- {\rm rk} (F) \bC]
\in  K^{-1}_{\bR/\bZ} (M)$, thus $( S^+(TZ)^*- S^-(TZ)^*)\otimes
\mathcal{F} \in   K^{-1}_{\bR/\bZ} (M)$. Set
\begin{align}\label{h3}
I(F)=  \sum_{i=0}^{n}(-1)^i\Big (H^i(Z,F|_Z), h^{H^i(Z,F|_Z)},
\nabla^{H^i(Z,F|_Z),e}, 0\Big ).
\end{align}
Let $\widehat{\eta} (\nabla ^F, h^F) \in \Omega ^{\rm odd}
(B)/{\rm Im} d$
 be the eta form of Bismut-Cheeger
\begin{align}\label{h4}
\widehat{\eta} (\nabla ^F, h^F)= \left({2\pi \sqrt{-1}}\right)^{-\frac{1}{2}}
\int_0^{\infty}\varphi  \tr_s\left[
\frac{\partial  C_{u}}{\partial u}
\exp (-  C_{u}^2)\right] du.
\end{align}
\begin{thm}
\begin{multline}\label{h5}
\ind_{top}(( S^+(TZ)^*- S^-(TZ)^*)\otimes
\mathcal{F} )\\
 = I(F) -\widehat{\eta} (\nabla ^F, h^F)
- {\rm rk} (F)\left(I(\bC)-\widehat{\eta} (\nabla ^\bC, h^\bC)\right).
\end{multline}
\end{thm}
\begin{proof} In \cite{L}, Lott used spin$^c$ Dirac operator to define
 $\ind_{an}$, especially, the spin$^c$ Dirac operator
(twisted by $( S^+(TZ)^*- S^-(TZ)^*)\otimes F$), $D^{Z,c}$ is
\begin{align}\label{h6}
D^{Z,c}= \sum_j c(e_j) \nabla ^{TZ\otimes F,e}_{e_j}
=D^Z +\frac{1}{2} \sum_j \widehat{c}(e_j) \omega(F,h^F)(e_j).
\end{align}
This operator $D^{Z,c}$ has not constant dimensional kernel, thus
we need to choose smooth finite dimensional sub-bundles $F_{\pm}$
of $E_{\pm}$ ($E_+= E^{\rm even}, E_- = E^{\rm odd}$) and
complementary subbundles $G_{\pm}$ such that $D^{Z,c}$ are
diagonal with respect to the decompositions $E_{\pm}=
F_{\pm}\oplus G_{\pm}$ and $D^{Z,c}$ restricted to $G_{\pm}$ is
invertible (cf. \cite[Definition 14]{L}). It seems that it is hard
to compare directly the right hand side of (\ref{h5}) to
$\ind_{an}(( S^+(TZ)^*- S^-(TZ)^*)\otimes \mathcal{F} )$ in
\cite[Definition 14]{L}. But by the arguments in \cite[Proposition
6 and Corollary 1]{L}, we will get (\ref{h5}) if we can prove the
following identity for any odd dimensional compact spin$^c$
manifold $B$,
\begin{align}\label{h7}
\overline{\eta}_M(( S^+(TZ)^*- S^-(TZ)^*)\otimes F) =
\overline{\eta}_B (I(F)),
\end{align}
where $\overline{\eta}$ is the reduced eta invariant of the
spin$^c$ Dirac operator twisted by the corresponding bundles as in
\cite[Definition 11]{L}.

Let $g^{TB}$ be a metric on $TB$ and let $\nabla ^{TB}$ be the Levi-Civita
connection on $(TB, g^{TB})$.  Let $\nabla ^{S(TB)}$ be the connection
on the spinor bundle $S(TB)$ induced by  $\nabla ^{TB}$
and the connection on the line budle defining the spin$^c$ structure.
 Let $\nabla ^e$ be the connection on
$\Lambda (T^* Z) \otimes \pi ^*(S(TB)) \otimes F$ induced by
$\nabla ^{\Lambda (T^* Z)}$, $\nabla ^{ \pi ^*(S(TB))}$, $\nabla ^{F,e}$. Set
\begin{align}\label{h8}
D^H= \sum_{\alpha} c(f_\alpha) (\nabla ^e_{f_\alpha}
 +\frac{1}{2} k(f_{\alpha})).
\end{align}
Then by proceeding as in \cite[(4.26)]{BC}, we get
\begin{align}\label{h9}
D^{M,c}= D^H + D^{Z,c} -\frac{c(T)}{4}.
\end{align}
and $\overline{\eta}_M(( S^+(TZ)^*- S^-(TZ)^*)\otimes F)$ is the reduced eta invariant of the operator $D^{M,c}$. Let
\begin{align}\label{h10}
D^{'M}= D^H + D^{Z} -\frac{c(T)}{4}.
\end{align}
Then $D^{'M}= D^{M,c}-\frac{1}{2} \sum_j \widehat{c}(e_j)
\omega(F,h^F)(e_j)$ and $\frac{1}{2} \sum_j \widehat{c}(e_j)
\omega(F,h^F)(e_j)$ 
anti-commutes  with
 $c(X),\ X\in TM$. Using the variation formula for  eta
invariants (cf. \cite{APS2} and \cite{BF}) and the local index
techniques as in \cite[Theorem 2.7]{BC}, we know (cf.
\cite[Proposition 3.5]{B2} and \cite[(3.5)]{Z})
\begin{align}\label{h11}
\overline{\eta}_M(( S^+(TZ)^*- S^-(TZ)^*)\otimes F) =
\overline{\eta}(D^{M,c}) = \overline{\eta}(D^{'M}) \quad \mod \bZ.
\end{align}
From (\ref{h11}), we can use the adiabatic limit argument as in
\cite[Proposition 6]{L} to get (\ref{h7}).
\end{proof}

By \cite[Theorem 3.7]{B2} and \cite[Proposition 2.3]{Z} (cf. also
Lemma \ref{t3.11}),
\begin{align}\label{h12}
\widehat{\eta} (\nabla ^F, h^F)=0.
\end{align}
Thus from (\ref{h4}) and (\ref{h5}), we know
\begin{align}\label{h13}
\ind_{top}(( S^+(TZ)^*- S^-(TZ)^*)\otimes
\mathcal{F} )
=I(F) - {\rm rk} (F) I(\bC) .
\end{align}
When we apply (\ref{h1})  to (\ref{h13}), we get
again (\ref{f6}).

Thus, (\ref{h13})
represents a refinement of (\ref{f6}) in $K^{-1}_{\bR/\bZ} (B)$.
We leave the interested reader to extend this to the case where no
spin$^c$ assumption on $TZ$ is required.


\subsection{ \normalsize The $\eta$  and torsion forms}
\label{s3.20}

$\quad$ In the rest of this section, the supertrace $\tr_s$ on $E$
is defined by the $\bZ_2$-grading induced from $(-1)^{N_Z}$. Let
\begin{align}\label{0d1}
d(H(Z,F)) =\sum_{i=0}^n (-1)^i i \dim H^i(Z,F).
\end{align}
Let $\widehat{\Lambda (T^*Z)}$ be another copy of $\Lambda (T^*Z)$.
For $\omega \in \Lambda (T^*Z)$, we denote by
$\widehat{\omega}\in \widehat{\Lambda (T^*Z)}$ the copy of $\omega$.
Then the Berezin integral $\int^B : \Lambda (T^*M)\widehat{\otimes}
\widehat{\Lambda (T^*Z)}\to \Lambda (T^*M)\otimes o(TZ)$ is defined by
$\int^B \gamma \wedge \widehat{\omega} \to  \gamma \,
i_{\widehat{e_1}}\cdots i_{\widehat{e_n}}\widehat{\omega} $
for $\gamma\in \Lambda (T^*M)$. We define also the fiber-wise
integral $\int_Z$ by: for $\gamma \in C^\infty(B, \Lambda (T^*B))$,
$\delta  \in C^\infty(M,\Lambda (T^*Z)\otimes o(TZ))$,
\begin{align}\label{0d2}
\int_Z (\pi ^* \gamma) \wedge \delta = \gamma \int_Z \delta.
\end{align}
Set
\begin{align}\label{0d3}
&\dot{R}{}^{TZ} =\sum_{i,j=1}^n \left\langle R^{TZ} e_i,e_j\right\rangle
\widehat{e^i}\wedge \widehat{e^j} \in \Lambda ^2(T^*M)\widehat{\otimes}
\widehat{\Lambda ^2 (T^*Z)},\\
&e(TZ,\nabla ^{TZ})= (-1)^{\frac{n(n+1)}{2}} \pi ^{\frac{n}{2}}
\int^B \exp\Big(-\frac{1}{2}\dot{R}{}^{TZ}\Big),\nonumber\\
&a_{-1}= (-1)^{\frac{n(n+1)}{2}} \pi ^{\frac{n}{2}}\rk(F)\, \int_Z
\int^B \frac{1}{2} \sum_{i=1}^n e^i\wedge \widehat{e^i}
\exp\Big(-\frac{1}{2}\dot{R}{}^{TZ}\Big). \nonumber
\end{align}
The form $e(TZ,\nabla ^{TZ})$ is the Chern-Weil representative of
the Euler class of $TZ$.  And $a_{-1}$ is a function on $B$ and is
$0$ if $n$ is even.

For any $u>0$, let $\psi_u: \Lambda (T^*B)\to \Lambda (T^*B)$ be
defined by that
 for $\gamma \in \Lambda (T^*B)$, $\psi_u\gamma =u^{-\deg \gamma/2} \gamma$.
Then by (3.14),
\begin{align}\label{1d4}
\tr_s \left[N_Z\exp((1+r^2)D_u^2)\right] = \psi_u \tr_s
\left[N_Z\exp((1+r^2)u D_1^2)\right].
\end{align}
By standard results on heat kernels, we know that $\tr_s
\left[N_Z\exp((1+r^2)u D_1^2)\right]$ has an asymptotic expansion
in $u$ as $u\to 0^+$, which only contains integral powers of $u$
if $n=\dim Z$ is even, and only contains half-integral powers of
$u$ if $n$ is odd. Since $\tr_s [N_Z\exp((1+r^2)D_u^2)]$ is an
even form on $B$, by (\ref{1d4}) we see that the same happens to
it.

 On the other hand, as in \cite[(11.1)]{BZ}, we have
\begin{align}\label{0d4}
N_Z=\frac{1}{2} \sum_{i=1}^n c(e_i) \widehat{c}(e_i)+\frac{n}{2}.
\end{align}
By \cite[Theorem 3.15]{BL} and (\ref{0d4}), as in \cite[Theorem
3.21]{BL}, \cite[Theorem 7.10]{BZ}, we have for $r\in \bR$, as
$u\to 0^+$,
\begin{align}\label{0d6}
\tr_s \left[N_Z\exp((1+r^2)D_u^2)\right] = \left
\{\begin{array}{l} \dip{\frac{n}{2} \chi (Z)\rk (F) +O(u) \quad
{\rm if}\ \ n\ \mbox{\rm is even},}\\ \\
\dip{\frac{a_{-1}}{\sqrt{(1+r^2)u}} + O(\sqrt{u})\quad {\rm if}\
\ n\ \mbox{\rm is odd}. }\end{array}\right.
\end{align}
While as $u\to  +\infty$,
\begin{align}\label{0d5}
\tr_s \left[N_Z\exp((1+r^2)D_u^2)\right] = d(H(Z,F))
+O\left(\frac{1}{\sqrt{u}}\right).
\end{align}
From (\ref{c21}) and (\ref{0d6}), we know that as $u\to 0^+$,
\begin{align}\label{1d5}
\tr_s\left[ \left(\frac{\partial}{\partial u} \left(C_{u}+
\sqrt{-1}r D_{u}\right) \right) \exp \left(-  (C_{u}+\sqrt{-1} r
D_{u})^2\right)\right]= \frac{\sqrt{-1}r
da_{-1}}{2\sqrt{1+r^2}u^{3/2}}
 + O\left(\frac{1}{\sqrt{u}}\right).
\end{align}

The following definition is closely related to \cite[Definition 3.22]{BL}.

\begin{defn}\label{t4.3}  For any $r\in{\bf R}$, put
\begin{multline}\label{0d7}
I_r=-r  
\varphi
\int_0^{+\infty}
\left(\tr_s\left[N_Ze^{-(C_u+\sqrt{-1}rD_u)^2}\right]- d(H(Z,F|_Z)) \right.\\
\left. -\frac{a_{-1}}{\sqrt{(1+r^2)u}}- \left(\frac{n}{2} \chi
(Z)\rk (F) -d(H(Z,F|_Z))\right)e^{-(1+r^2)u/4}\right){du\over 2u}.
\end{multline}
\end{defn}
Let
 $\widehat{\eta}_r$ be the $\eta$-form of Bismut-Cheeger
\cite[Definition 4.33]{BC} defined by
\begin{multline}\label{d1}
\widehat{\eta}_r=  \left({2\pi \sqrt{-1}}\right)^{-\frac{1}{2}}
\int_0^{\infty}\varphi \Big \{  \tr_s\left[
\left(\frac{\partial}{\partial u} \left(C_{u}+ \sqrt{-1}r
D_{u}\right) \right) \exp \left(-  (C_{u}+\sqrt{-1} r
D_{u})^2\right)\right]\\
 -\frac{\sqrt{-1}r da_{-1}}{2\sqrt{1+r^2}u^{3/2}} \Big \} du.
\end{multline}

\begin{Rem}
The extra term involving $da_{-1}$ in the right hand side of
(\ref{d1}) shows that this $\widehat{\eta}_r$ form is slightly
different from what in \cite{BC}.
\end{Rem}

\begin{thm}\label{t4.1} For any $r\in {\bf R}$, the $\eta$ form
$\widehat{\eta}_r$ is exact and is $0$ at $r=0$. Moreover,  the
 following transgression formula holds,
\begin{align}\label{0d8}
\widehat{\eta}_r=-{1\over 2\pi}dI_r.
\end{align}
\end{thm}
\begin{proof}
Theorem \ref{t4.1} is a direct consequence of  Lemma \ref{t3.11}
and (\ref{0d5}), (\ref{0d6}).
\end{proof}

Let $T_f(T^HM,g^{TZ}, h^F)$ be the torsion form constructed in
the spirit of  \cite[Definition 3.21]{BL}
associated to the odd holomorphic function $f(z)$ such that $f'(z)=e^{z^2}$,
that is
\begin{multline}\label{0d9}
T_f(T^HM,g^{TZ}, h^F)=-
\varphi \int_0^{+\infty}
\left(\tr_s\left[N_Ze^{D_u^2}\right]\right. - d(H(Z,F|_Z))\\
\left. - \frac{a_{-1}}{\sqrt{u}}
 -\left(\frac{n}{2} \chi (Z)\rk (F) - d(H(Z,F|_Z))\right)e^{-u/4}\right){du\over
 2u}.
\end{multline}

\begin{thm}\label{t4.5}   The following identity holds,
\begin{align}\label{0d10}
\left. {\partial I_r\over \partial r}\right|_{r=0}=
T_f(T^HM,g^{TZ}, h^F).
\end{align}
In particular,
\begin{align}\label{0d11}
\left. {\partial \widehat{\eta}_r\over \partial
r}\right|_{r=0}=-{1\over 2\pi}dT_f(T^HM,g^{TZ}, h^F).
\end{align}
\end{thm}

\begin{proof} Formula (\ref{0d10}) follows from (\ref{0d7}), (\ref{0d9}).
 Formula (\ref{0d11}) follows from (\ref{0d8}) and (\ref{0d10}).
\end{proof}


\begin{thm}\label{t4.6} For any $r\in{\bf R}$, the following identity holds,
\begin{multline}\label{0d12}
\widehat{\eta}_r={r\over 2\pi}\sum_{j=0}^{+\infty}{(1+r^2)^j\over
j!({2j+1})}
\sum_{i=1}^n{(-1)^i}c_{2j+1}(H^i(Z,F|_Z),h^{H^i(Z,F|_Z)})\\
-{r\over 2\pi}\sum_{j=0}^{+\infty}{(1+r^2)^j\over j! (2j+1)}
\int_Z e(TZ,\nabla ^{TZ})c_{2j+1}(F,h^F).
\end{multline}
In particular,
\begin{multline}\label{0d13}
\left. {\partial\widehat{\eta}_r\over \partial  r}\right|_{r=0}=
{1\over 2\pi}\sum_{j=0}^{+\infty} {1\over
j!{(2j+1)}}\sum_{i=1}^n{(-1)^i}c_{2j+1}(H^i(Z,F|_Z),h^{H^i(Z,F|_Z)})\\
- {1\over 2\pi}\sum_{j=0}^{+\infty} {1\over
j!(2j+1)} \int_Z e(TZ,\nabla ^{TZ})c_{2j+1}(F,h^F).
\end{multline}
\end{thm}
\begin{proof}
In fact,  by Lemma \ref{t3.11} and (\ref{d1}), we have
\begin{align}\label{0d14}
\widehat{\eta}_r =\sqrt{-1}r \left({2\pi
\sqrt{-1}}\right)^{-\frac{1}{2}}  \int_0^{\infty}\varphi \left\{
\tr_s\left[\frac{\partial D_{u}}{\partial u} \exp \left( (1+r^2)
D_{u}^2\right)\right] -\frac{ da_{-1}}{2\sqrt{1+r^2}u^{3/2}}
\right\} du.
\end{align}

From (\ref{a8}) and (\ref{0d14}), one deduces that
\begin{align}\label{0d15}
\left.\widehat{\eta}_r = \frac{r}{\sqrt{1+r^2}}
\psi_{(1+r^2)^{-1}} {\partial\widehat{\eta}_r\over \partial  r}\right|_{r=0}.
\end{align}
By (\ref{0d15}), we  need only to prove (\ref{0d13}).

\begin{lemma}\label{t4.001} The following identity holds,
\begin{align}\label{0d16}
\tr_s\left[\frac{\partial D_u}{\partial u} \exp(D_u^2)\right]
-\frac{da_{-1}}{2u^{3/2}} = \frac{\partial}{\partial u}\int_0^1
\tr_s\left[D_u\exp(s^2 D_u^2)\right] ds.
\end{align}
\end{lemma}

\begin{proof} By (\ref{a8}), (\ref{c21}) and (\ref{0d6}), one has
\begin{multline}\label{0d160}
\lim_{s\to 0^+}\tr_s\left[ \frac{\partial D_u}{\partial
u}s\exp(s^2 D_u^2)\right ] =\lim_{s\to 0^+} {s\over 2u}d
\psi_{s^{-2}}\tr_s\left[ N_Z\exp( D_{s^2u}^2)\right
]\\
=\frac{da_{-1}}{2u^{3/2}}.
\end{multline}
Thus
\begin{multline}\label{0d17}
\frac{\partial}{\partial u}\int_0^1 \tr_s[D_u\exp(s^2 D_u^2)] ds\\
= \int_0^1 \tr_s\left[\frac{\partial D_u}{\partial u}\exp(s^2 D_u^2)\right] ds
+ \int_0^1 \tr_s\left[D_u \left[s^2 D_u, \frac{\partial D_u}{\partial u}\right]\exp(s^2 D_u^2)\right] ds\\
=\int_0^1 \tr_s\left[\frac{\partial D_u}{\partial u} (1+ 2s^2 D_u^2)\exp(s^2 D_u^2)\right] ds\\
=  \tr_s\left[ \frac{\partial D_u}{\partial u}\int_0^1
\frac{\partial}{\partial s}(s\exp(s^2 D_u^2))ds\right]\\
 =\tr_s\left[\frac{\partial D_u}{\partial u} \exp(D_u^2)\right]-
\frac{da_{-1}}{2u^{3/2}},
\end{multline}
which is exactly (\ref{0d16}).
\end{proof}

Now by (\ref{a8}) again, one has
\begin{align}\label{0d18}
 \tr_s\left[D_u\exp(s^2 D_u^2)\right]
=  {1\over s}\psi_{s^{-2}}
\tr_s\left[D_{s^2u}\exp(D_{s^2u}^2)\right] .
\end{align}
Now from \cite[Theorem 3.16]{BL}, we know that
\begin{multline}\label{0d19}
(2 \pi \sqrt{-1})^{\frac{1}{2}} \varphi \tr_s\left[D_{u}\exp(D_{u}^2)\right] \\
=\left\{ \begin{array}{l}
\dip{\int_Z e(TZ,\nabla ^{TZ})
\sum_{j=0}^\infty \frac{1}{j!}c_{2j+1}(F,h^F)
+ O(\sqrt{u})\quad  {\rm as} \ u\to 0,}\\
\dip{  \sum_{j=0}^\infty \frac{1}{j!}\sum_{i=0}^n(-1)^ic_{2j+1}(H^i(Z,F|_Z),h^{H^i(Z,F|_Z)})
+ O\left(\frac{1}{\sqrt{u}}\right) }\\
\dip{\hspace*{50mm} {\rm as} \ u\to +\infty.}
\end{array}\right.\end{multline}

Since $\tr_s\left[D_u\exp(s^2 D_u^2)\right]$ is an odd form on
$B$, by (\ref{0d19}), one sees that $|{1\over s}\psi_{s^{-2}}
\tr_s[D_{s^2u}\exp(D_{s^2u}^2)]|$ has a fixed uniform upper bound
for
 $s\in(0,1], \ u\in (0,+\infty]$.

  Thus, from (\ref{0d14}), (\ref{0d16}),
(\ref{0d18}), (\ref{0d19}) and the dominated convergence property,
we get
\begin{multline}\label{0d20}
\left. 2\pi {\partial\widehat{\eta}_r\over \partial  r}\right|_{r=0}=
(2 \pi \sqrt{-1})^{1/2} \varphi
\left[\lim_{u\to +\infty} \int_0^1{1\over s} \psi_{s^{-2}} \tr_s\left[D_{s^2u}\exp(D_{s^2u}^2)\right] ds\right.\\
\left.-\lim_{u\to 0^+} \int_0^1 {1\over s} \psi_{s^{-2}}
\tr_s\left[D_{s^2u}\exp(D_{s^2u}^2)\right] ds\right]\\
=      \sum_{j=0}^\infty
\frac{1}{j!}\int_0^1s^{2j}ds\sum_{i=0}^n(-1)^ic_{2j+1}(H^i(Z,F|_Z),h^{H^i(Z,F|_Z)})\\
-  \int_Z e(TZ,\nabla ^{TZ}) \sum_{j=0}^\infty
\frac{1}{j!}c_{2j+1}(F,h^F)\int_0^1s^{2j}ds\\
=\sum_{j=0}^{+\infty} {1\over
j!{(2j+1)}}\sum_{i=1}^n{(-1)^i}c_{2j+1}(H^i(Z,F|_Z),h^{H^i(Z,F|_Z)})\\
- \sum_{j=0}^{+\infty} {1\over j!(2j+1)} \int_Z e(TZ,\nabla
^{TZ})c_{2j+1}(F,h^F),
\end{multline}
which is equivalent to (\ref{0d13}).
\end{proof}

Combining (\ref{0d11}) and (\ref{0d13}), one gets the following
transgression formula of Bismut-Lott type.

\begin{cor}\label{t5} The following identity holds,
\begin{multline}\label{0d33}
dT_f(T^HM,g^{TZ}, h^F)=\sum_{j=0}^{+\infty} {1\over j!(2j+1)}
\int_Z e(TZ,\nabla ^{TZ})c_{2j+1}(F,h^F)\\
-\sum_{j=0}^{+\infty} {1\over
j!{(2j+1)}}\sum_{i=1}^n{(-1)^i}c_{2j+1}(H^i(Z,F|_Z),h^{H^i(Z,F|_Z)}).
\end{multline}
\end{cor}

Alternately, (\ref{0d33}) is also a consequence of the following
more precise relation between $T_f(T^HM,g^{TZ}, h^F)$ and the
Bismut-Lott torsion form ${\cal T}(T^HM,g^{TZ}, h^F)$  defined in
\cite[Definition 3.22]{BL}.
\begin{thm} The following identity holds in $\Omega^*(B)$,
\begin{align}\label{0d21}
{\cal T}(T^HM,g^{TZ}, h^F) = (1+N_B) T_f(T^HM,g^{TZ}, h^F).
\end{align}
\end{thm}
\begin{proof} Recall that the Bismut-Lott torsion form ${\cal T}(T^HM,g^{TZ}, h^F)$
is defined by
\begin{multline}\label{0d000}
{\cal T}(T^HM,g^{TZ}, h^F)=-\varphi
\int_0^{+\infty} \left(\tr_s\left[N_Z(1+2D_u^2)e^{D_u^2}\right]\right.\\
- d(H(Z,F|_Z)) \left. - \left(
  \frac{n}{2} \chi (Z)\rk (F) - d(H(Z,F|_Z))\right)\left(1-{u\over 2}\right)e^{-u/4}\right){du\over  2u}.
\end{multline}

A direct computation shows that the 0-form component of
$T_f(T^HM,g^{TZ}, h^F)$ is exactly the half of the Ray-Singer
analytic torsion defined in \cite{RS} and \cite{BZ}. Thus,  the
$0$-form component of (\ref{0d21}) is a consequence of
\cite[Theorem 3.29]{BL}.

On the other hand, for $i>0$, we denote by a superscript $[i]$ the
$i$-form component of the corresponding forms. Then by (\ref{a8}),
one has
\begin{align}\label{1d21}
\left\{  \tr_s\left[ N_Z D_u^2 \exp \left(D_u^2\right)\right]
\right\} ^{[i]} = u^{-i/2}\left\{  \tr_s\left[ N_Z u D_1^2 \exp
\left(u D_1^2\right)\right] \right\} ^{[i]}.
\end{align}
Thus, one deduces that
\begin{multline}\label{1d22}
\int_0^{\infty} \left\{  \tr_s\left[ N_Z D_u^2 \exp
\left(D_u^2\right)\right] \right\} ^{[i]} \frac{du}{u}
=\int_0^{\infty} u^{-\frac{i}{2}} \left\{ \tr_s\left[ N_Z  D_1^2
\exp \left(u
D_1^2\right)\right] \right\} ^{[i]} du\\
=\int_0^{\infty} u^{-\frac{i}{2}}\frac{\partial}{\partial u}
\left\{ \tr_s\left[ N_Z  \exp \left(u D_1^2\right)\right] \right\}
^{[i]} du  \\
=\int_0^{\infty} i\, u^{-\frac{i}{2}} \left\{ \tr_s\left[ N_Z \exp
\left(u D_1^2\right)\right] \right\} ^{[i]} \frac{du}{2u},
\end{multline}
where in the last equality we have used the facts that
\begin{multline}\label{1d220}
\lim_{u\rightarrow 0^+}u^{-\frac{i}{2}}\left\{ \tr_s\left[ N_Z
\exp \left(u D_1^2\right)\right] \right\}
^{[i]}\\
=\lim_{u\rightarrow 0^+}\left\{ \tr_s\left[ N_Z  \exp \left(
D_u^2\right)\right] \right\} ^{[i]} =0,
\end{multline}
and
\begin{multline}\label{1d221}
\lim_{u\rightarrow +\infty}u^{-\frac{i}{2}}\left\{ \tr_s\left[ N_Z
\exp \left(u D_1^2\right)\right] \right\}
^{[i]}\\
=\lim_{u\rightarrow +\infty}\left\{ \tr_s\left[ N_Z  \exp \left(
D_u^2\right)\right] \right\} ^{[i]} =0,
\end{multline}
which are the consequences of (\ref{0d6}) and (\ref{0d5}).

 From
(\ref{1d22}), we get (\ref{0d21}).
\end{proof}

\begin{Rem}
From (\ref{c21}),  as in (\ref{d10}),  (\ref{d11}), we get
\begin{multline}\label{d12}
\frac{\partial}{\partial r} \left [\frac{r}{2u}d \tr_s
\left[N_Z \exp \left(-  (C_{u}+\sqrt{-1} r D_{u})^2\right)\right]\right]\\
- \frac{\partial}{\partial u}\tr_s\left[D_{u}\exp \left(-  (C_{u}+\sqrt{-1} r D_{u})^2\right)\right]\\
= \frac{-1}{2u}d \tr_s\left[ 2 N_Z D_{u}^2 \exp \left(-  (C_{u}+\sqrt{-1} r D_{u})^2\right)\right].
\end{multline}
Especially, when we restrict ourselves  to $r=0$, from
(\ref{d12}), we get
\begin{align}\label{d13}
 \frac{\partial}{\partial u}
\tr_s\left[D_{u}\exp \left(- C_{u}^2\right)\right]
=\frac{1}{2u}  d \tr_s\left[ N_Z (1+2D_u^2) \exp \left(- C_u^2\right)\right].
\end{align}
This is exactly \cite[Theorem 3.20]{BL}. It seems interesting
that here we obtain it purely
 through the consideration of $\eta$ forms.
\end{Rem}

\end{document}